\documentclass[11pt]{article}
\usepackage{jmlr2e-mod}
\pdfoutput=1

\usepackage{color}
\usepackage{amsmath,amsfonts,amsbsy,amssymb}
\usepackage{booktabs}
\usepackage{caption,subcaption}
\usepackage{colortbl}
\usepackage{float}
\usepackage{fixmath} 
\usepackage{graphicx}
\usepackage[htt]{hyphenat} 
\usepackage{ifthen}
\usepackage{ltablex}
\usepackage{mathtools}
\usepackage[mathscr]{eucal} 
\usepackage{multirow}
\usepackage{pgffor}
\usepackage{siunitx}
\usepackage{stmaryrd} 
\usepackage[table]{xcolor}
\usepackage{tabularx}
\usepackage{url}
\usepackage{wrapfig}
\graphicspath{{./}}

\title{Parameter Sensitivity Analysis of the SparTen High Performance 
	Sparse Tensor Decomposition Software: Extended Analysis}

\author{\name Jeremy M. Myers \email 
jermyer@sandia.gov, jmmyers01@email.wm.edu 
\\
	\name Daniel M. Dunlavy \email dmdunla@sandia.gov \\
	\name Keita Teranishi \email knteran@sandia.gov \\
	\name D. S. Hollman \email dshollm@sandia.gov \\
	\addr Sandia National Laboratories, Albuquerque, NM 87123, USA \\
	\addr College of William and Mary, Williamsburg, VA 23185, USA}

\begin{document}
		
\maketitle

\renewcommand*{\thefootnote}{(\fnsymbol{footnote})}
\renewcommand*{\thefootnote}{\arabic{footnote}.}
	
\begin{abstract}
Tensor decomposition models play an increasingly important role in modern 
data science applications. One problem of particular interest is fitting a 
low-rank Canonical Polyadic (CP) tensor decomposition 
model when the tensor has sparse structure and the tensor elements are 
nonnegative count data. SparTen is a high-performance C++ library which 
computes a low-rank decomposition using different solvers: a first-order 
quasi-Newton or a second-order damped Newton method, along with the appropriate 
choice of runtime parameters. Since default parameters in SparTen are tuned to 
experimental results in prior published work on a single real-world dataset conducted 
using MATLAB implementations of these methods, it remains unclear if the 
parameter defaults in SparTen are appropriate for general tensor data. 
Furthermore, it is unknown how sensitive algorithm convergence is to changes in 
the input parameter values. This report addresses these unresolved issues with 
large-scale experimentation on three benchmark tensor data sets. Experiments 
were conducted on several different CPU architectures and replicated with many 
initial states to establish generalized profiles of algorithm convergence behavior. 
\end{abstract}

\begin{keywords}
tensor decomposition, Poisson factorization, Kokkos, Newton optimization
\end{keywords}

\begin{figure}[b!]
\centering
\includegraphics[width=\textwidth]{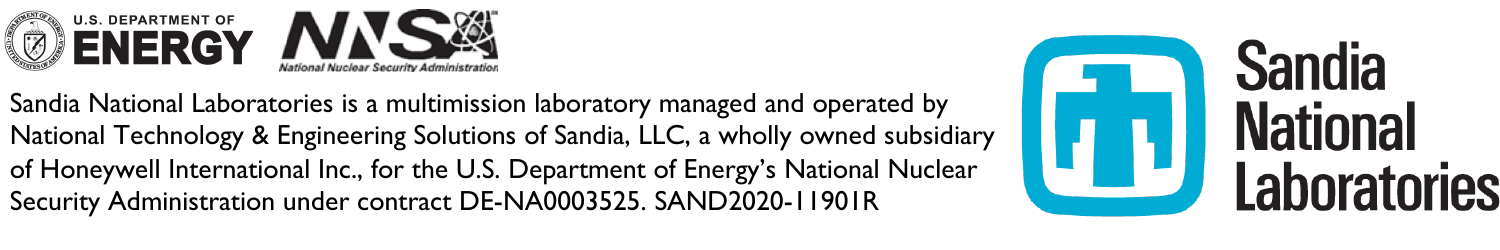}
\end{figure}

\clearpage
\section*{Acknowledgment}
We would like to thank Richard Barrett for assistance utilizing computing 
resources at Sandia National Laboratories and Rich Lehoucq for comments of 
support.

\tableofcontents
\listoffigures
\listoftables
\clearpage

\section{Introduction}
\label{sec:intro}
The Canonical Polyadic (CP) tensor decomposition model has garnered attention 
as a tool for extracting useful information from high dimensional data across a 
wide range of 
applications~\cite{doi:10.1137/07070111X,7038247,doi:10.1002/sapm192761164,CarrollChang70,harshman70}.

Recently, Hansen {\it et al.} developed two highly-parallelizable Newton-based 
methods for low-rank tensor factorizations on Poisson count data 
in~\cite{HaPlKo15}, one a first-order quasi-Newton method (PQNR) and another a 
second-order damped Newton method (PDNR). They were first implemented in MATLAB 
Tensor Toolbox~\cite{TTB_Software} as the function \verb|cp_apr|, referring to 
this approach as computing a CP decomposition using Alternating Poisson 
Regression (i.e., CP-APR). These methods fit a reduced-rank CP model to count 
data, assuming a Poisson error distribution. PDNR and PQNR are implemented in 
SparTen,\footnote{SparTen is a portmanteau 
word derived from \emph{Sparse} and \emph{Tensor}. The SparTen code is available at 
\texttt{http://gitlab.com/tensors/sparten}.} a 
high-performance C++ library of CP-APR solvers for sparse tensors. SparTen 
improves on the MATLAB implementation to provide efficient execution for large, 
sparse tensor decompositions, exploiting the Kokkos hardware abstraction 
library~\cite{CarterEdwards20143202} to harness parallelism on diverse HPC 
platforms, including x86-multicore, ARM, and GPU computer architectures.

SparTen contains many algorithmic parameters for controlling the optimization 
subroutines comprising PDNR and PQNR. To date, only anecdotal 
evidence exists for how best to tune the algorithms. Parameter defaults in 
SparTen were chosen according to previous results using the MATLAB 
implementations described by Hansen {\it et al.}~\cite{HaPlKo15}. However, 
their analysis was limited to a 
single real-world dataset, and thus may not be optimal for computing 
decompositions of more general tensor data. Furthermore, it is unknown how the 
initial guess to a solution affects convergence, since SparTen methods may 
converge slowly---or worse, stagnate---on real data if the initial state is far 
from a solution. And, lastly, the average impact of input parameters on algorithm 
convergence is unclear.

To address these unknowns, we present the results of numerical experiments to 
assess the sensitivity of software parameters on algorithm convergence for a 
range of values with benchmark tensor problems. Every experiment was replicated 
with 30 randomly chosen initial guesses on three diverse computer 
architectures to aid statistical interpretation. With our results, we (1) 
provide new results that offer a realistic picture of algorithm convergence 
under reasonable resource constraints, (2) establish practical bounds on 
parameters such that, if set at or beyond these values, convergence is 
unlikely, and (3) identify areas of performance degradation and convergence 
toward qualitatively different results owing to parameter sensitivities.

We limited our study to multicore CPU architectures only, using 
OpenMP~\cite{dagum1998openmp} to manage the parallel computations across 
threads/cores. Although SparTen, through Kokkos, can leverage other execution 
backends---e.g., NVIDIA's CUDA framework for GPU computation---we focus solely 
on diversity in CPU architectures in this work.

This paper is structured as follows. Section~\ref{sec:background} summarizes 
basic tensor notation and details. Section~\ref{sec:methods} describes
the hardware environment, test data, and experimental design of the 
sensitivity analysis. Section~\ref{sec:results} provides detailed results of 
the sensitivity analyses. Section~\ref{sec:conc} offers concluding remarks and 
lays out future work.
\clearpage

\section{Background}
We briefly describe below the problem we are addressing in this report; for a 
detailed description of CP decomposition algorithms implemented in SparTen, 
refer to the descriptions in Hansen {\it et al.}~\cite{HaPlKo15}.

An $N$-way data tensor $\boldsymbol{\mathscr{X}}$ has dimension sizes $I_1 
\times I_2 \times \dots \times I_N$. We wish to fit a reduced-dimension tensor 
model, $\boldsymbol{\mathscr{M}}$, to $\boldsymbol{\mathscr{X}}$. The 
$R$-component Canonical Polyadic (CP) decomposition is given as follows: 
\begin{equation} 
\boldsymbol{\mathscr{X}} \approx \boldsymbol{\mathscr{M}} = \llbracket 
\mathbold{\lambda}; 
\mathbf{A}^{(1)},\ldots,\mathbf{A}^{(N)}\rrbracket = \sum_{r=1}^R \lambda_r 
\mathbf{a}_r^{(1)} \circ \ldots \circ \mathbf{a}_r^{(N)},\label{eq:cpmodel}
\end{equation}
where $\mathbold{\lambda} = [\lambda_1,\ldots,\lambda_R]$ is a scaling vector, 
$\mathbf{a}_r^{(n)}$ represents 
the $r$-th 
column of the factor matrix $\mathbf{A}^{(n)}$ of size $I_n \times R$, and 
$\circ$ is the vector outer product. We refer to the operator $\llbracket \cdot 
\rrbracket$ 
as a Kruskal operator, and the tensor $\boldsymbol{\mathscr{M}}$, with its 
specific multilinear model form, as a Kruskal tensor in~\eqref{eq:cpmodel}. 
See \cite{doi:10.1137/07070111X} for more details regarding these definitions.

SparTen addresses the special case when the elements of 
$\boldsymbol{\mathscr{X}}$ are nonnegative counts. Assuming the entries in 
$\boldsymbol{\mathscr{X}}$ follow a Poisson distribution with multilinear 
parameters, the low-rank CP decomposition in~\eqref{eq:cpmodel} can be computed 
using the CP-APR methods, PDNR and PQNR, introduced by Hansen {\it et al.}~\cite{HaPlKo15}. 
\label{sec:background}
\clearpage

\section{Methods}
\label{sec:methods}
In this section, we describe the hardware platforms, data, 
and SparTen algorithm parameters used in our experiments.

\subsection{Hardware Platforms}
\label{sec:methods:hardware}
We used diverse computer architectures running Red Hat Enterprise Linux (RHEL) 
to perform our experiments, with 
hardware and compiler specifications detailed in Table~\ref{tab:methods:hw}. 
Intel 1--4 are production clusters with hundreds to thousands of nodes, whereas 
ARM and IBM clusters are advanced architecture research testbeds with tens of 
nodes each. The ARM and IBM testbeds have larger memory and 
support many more threads per node than do any cluster in Intel 1-4. We 
employed the maximum number of OpenMP threads available per node 
from each platform to maximize throughput and configured the maximum wall-clock 
limit as 12:00 hours for all experiments. All parallelism was solely across 
threads on a single node. We built and compiled SparTen to leverage OpenMP via 
Kokkos with the latest 
software build tools available on each cluster. The GNU compiler, \texttt{gcc}, 
was used, with \texttt{-O3} 
optimization and Kokkos architecture-specific flags enabled.
\begin{table}[!h]
\centering
\caption{Hardware characteristics and software environment of the clusters 
in this paper. {\it Threads} and {\it RAM (GB)} are per node.}
\label{tab:methods:hw}
\begin{tabular}{llrrrrrr}
\toprule
Platform & Processor & Nodes & CPUs & Threads & RAM (GB) & GCC \\ 
\midrule
ARM 	& ThunderX2 	& 44 	& 28 		& 256 	& 255 & 7.2.0 \\
IBM 	& Newell Power9 & 10 	& 20 		& 80 	& 319 & 7.2.0 \\
Intel 1 & Sandy Bridge 	& 1,848 & 29,568 	& 16 	& 64  & 8.2.1 \\
Intel 2 & Broadwell 	& 740 	& 26,640 	& 36 	& 128 & 8.2.1 \\
Intel 3 & Sandy Bridge 	& 201 	& 3,344 	& 16 	& 64  & 8.2.1 \\ 
Intel 4 & Sandy Bridge 	& 1,232 & 19,712 	& 16 	& 64  & 8.2.1 \\
\bottomrule
\end{tabular}
\end{table}

\subsection{Data}
\label{sec:methods:datasets}
We conducted experiments using sparse tensors of count data from the FROSTT 
collection~\cite{frosttdataset}. Specifically, we chose the following three 
datasets 
(summary statistics provided in Table~\ref{tab:methods:datasets}) to account 
for size, dimensionality, 
and density (i.e., the ratio of nonzero entries to the total number of elements 
in the tensor):
\begin{enumerate}
\item Chicago Crime Community is a 4th-order tensor 
of crime reports in the city of Chicago spanning nearly 17 years. The four 
modes represent {\it day} $\times$ {\it hour} $\times$ {\it community} $\times$ 
{\it crime-type} and the values are counts. 
\item LBNL-Network is a 5th-order tensor of anonymized 
internal network traffic at Lawrence Berkeley National Laboratory. The five 
modes represent {\it sender-IP} $\times$ {\it sender-port} $\times$ {\it 
destination-IP} $\times$ {\it destination-port} $\times$ {\it time} and the 
values are total packet length per timestep. 
\item NELL-2 is 3rd-order benchmark tensor that gives a snapshot 
of the NELL: Never-Ending Language Learning relational database. The three 
modes represent {\it entity} $\times$ {\it relation} $\times$ {\it entity} 
relationships.
\end{enumerate} 
Throughout the discussion below, we refer to the data using the short names 
listed in the table.
\begin{table}[!h]
\centering
\caption{Sparse tensor datasets from the FROSTT collection.}
\label{tab:methods:datasets}
\resizebox{\columnwidth}{!}{%
\begin{tabular}{lrrr}
\toprule
FROSTT Name  (short name) & Nonzeros & Dimensions & Density \\
\midrule
Chicago Crime Community ({\it chicago})& 5.3M & $(6186,24,77,32)$ & $1.5 
\times 
10^{-02}$	\\
LBNL-Network ({\it lbnl})& 1.7M & $(1605,4198,1631,4209,868131)$ & $4.2 
\times 
10^{-14}$ \\
NELL-2 ({\it nell})& 77M & $(12092,9184,28818)$ & $2.4 \times 10^{-05}$ 
\\
\bottomrule
\end{tabular}
}%
\end{table}

\subsection{Software Parameter Definitions \& Experimental Ranges}
\label{sec:methods:params:description}
PQNR and PDNR are composed of standard techniques in the numerical optimization 
literature. Specifically, for each tensor mode, the Newton optimization 
computes the gradient and Hessian matrix. Then, the inverse Hessian is 
approximated to compute a search direction and an Armijo backtracking line 
search is used to compute the Newton step. How the inverse Hessian is 
approximated differentiates PDNR and PQNR. PDNR shifts the eigenvalues by a 
damping factor $\mu$ to guarantee the Hessian matrix is semi-positive definite, 
and solves the resulting linear system exactly. PQNR approximates the inverse 
Hessian directly with a limited-memory BFGS (L-BFGS) approach, computed with a 
small number of update pairs. Since the algorithm parameters analyzed here are 
those presented in several equations and algorithms in~\cite{HaPlKo15}, we 
defer to that paper for specific details.

To support discussion later in Section~\ref{sec:results}, we  
group the algorithm parameters into the following five categories. We note that 
the stability parameters used to safeguard against numerical errors---e.g., 
offset tolerances to avoid divide-by-zero floating point errors---do not appear 
in the corresponding Matlab Tensor Toolbox method \verb|cp_apr|.

\paragraph{A. CP-APR}
\begin{itemize}
	\item \texttt{max\_outer\_iterations}: Maximum number of outer iterations 
	to perform (Algorithm 1, Steps 2-9 in~\cite{HaPlKo15}).
	\item \texttt{max\_inner\_iterations}: Maximum number of inner iterations 
	to perform ($K_{max}$ in Algorithms 3 and 4 in~\cite{HaPlKo15}).
\end{itemize}
\paragraph{B. Line search}
\begin{itemize}
	\item \texttt{max\_backtrack\_steps}: Maximum number of backtracking steps 
	in line search (maximum allowable value of $t$ used in Equation (17) 
	in~\cite{HaPlKo15}).
	\item \texttt{min\_variable\_nonzero\_tolerance}: Tolerance for nonzero 
	line search step length (smallest allowable value of $\beta$ in Equation 
	(17) in~\cite{HaPlKo15}).
	\item \texttt{step\_reduction\_factor}: Factor to reduce line search step 
	size between iterations ($\beta^{t+1}/\beta^{t}$ in Equation (17) 
	in~\cite{HaPlKo15}).
	\item \texttt{suff\_decrease\_tolerance}: Tolerance to ensure the next 
	iterate decreases the objective function ($\sigma$ in Equation (17) 
	in~\cite{HaPlKo15}).
\end{itemize}
\paragraph{C. Damped Newton (PDNR)}
\begin{itemize}
	\item \texttt{mu\_initial}: Initial value of damping parameter ($\mu_0$ in 
	Algorithm 3 in ~\cite{HaPlKo15}).
	\item \texttt{damping\_increase\_factor}: Scalar value to increase damping 
	parameter in next iterate (Equation (16) in~\cite{HaPlKo15}).
	\item \texttt{damping\_decrease\_factor}: Scalar value to decrease damping 
	parameter in next iterate (Equation (16) in~\cite{HaPlKo15}).
	\item \texttt{damping\_increase\_tolerance}: Tolerance to increase the 
	damping parameter (Equation (16) in~\cite{HaPlKo15}). If the search 
	direction increases the objective function and the ratio of actual 
	reduction and predicted reduction in objective function ($\rho$ in Equation 
	(15) in~\cite{HaPlKo15}) is less than 
	\texttt{damping\_increase\_tolerance}, the damping parameter $\mu_k$ will 
	be increased for the next iteration.
	\item \texttt{damping\_decrease\_tolerance}: Tolerance to decrease the 
	damping parameter (Equation (16) in~\cite{HaPlKo15}). Conversely, if the 
	search direction decreases the objective function and the ratio of actual 
	reduction to predicted reduction ($\rho$ in Equation (15) 
	in~\cite{HaPlKo15}) is greater than \texttt{damping\_decrease\_tolerance}, 
	the damping parameter $\mu_k$ will be decreased for the next iteration.
\end{itemize}
\paragraph{D. Quasi-Newton (PQNR)}
\begin{itemize}
	\item \texttt{size\_LBFGS}: 
	Number of recent limited-BFGS (L-BFGS)-update pairs to use in estimating 
	the current Hessian ($M$ in Equation (18) in~\cite{HaPlKo15}).
\end{itemize}
\paragraph{E. Numerical stability}
\begin{itemize}
	\item \texttt{eps\_div\_zero\_grad}: Safeguard against divide-by-zero in 
	gradient and Hessian calculations.
	\item \texttt{log\_zero\_safeguard}: Tolerance to avoid computing $\log(0)$ 
	in objective function calculations.
\end{itemize}

The default value in SparTen of each parameter described above 
and the experimental ranges tested
in these experiments are given in Table~\ref{tab:methods:params}. 
\begin{table*}[ht!]
\centering
\caption{SparTen software parameter descriptions and values used in our 
experiments.}
\label{tab:methods:params}
\resizebox{\columnwidth}{!}{%
\begin{tabular}{lll}
\toprule
Parameter & Default & Values Used in Experiments \\ 
\midrule
\texttt{max\_outer\_iterations}$^1$
	& $10^5$ 
	& 1, 2, 4, 8, 16, 32, 64, 128, 256, 512 \\
\texttt{max\_inner\_iterations}$^1$
	& 20  
	& 20, 40, 80, 160 \\
\texttt{max\_backtrack\_steps}$^2$
	& 10 		
	& $1^\dagger, 2, 4, 8, 10, 12^\dagger, 16$ \\
\texttt{min\_variable\_nonzero\_tolerance}$^2$
	& $10^{-7}$ 
	& $10^{-1\dagger}, 10^{-3}, 10^{-7}, 10^{-15\dagger}$ \\
\texttt{step\_reduction\_factor}$^2$ 			
	& 0.5 
	& $0.1, 0.3^\dagger, 0.5, 0.7^\dagger, 0.9$ \\
\texttt{suff\_decrease\_tolerance}$^2$			
	& $10^{-4}$ 
	& $10^{-2}, 10^{-4}, 10^{-8\dagger}, 10^{-12\dagger}$ \\
\texttt{mu\_initial}$^3$			
	& $10^{-5}$	
	& $10^{-2}, 10^{-5},10^{-8}$ \\ 
\texttt{damping\_increase\_factor}$^3$
	& 3.5		
	& $1.5, 2.5^\dagger, 3.5, 4.5^\dagger, 5.5 $\\
\texttt{damping\_decrease\_factor}$^3$
	& $2/7$		
	& $0.1, 2/7, 0.3^\dagger, 0.5, 0.7^\dagger, 0.9$ \\
\texttt{damping\_increase\_tolerance}$^3$
	& 0.25 		
	& $0.1, 0.25, 0.495$ \\
\texttt{damping\_decrease\_tolerance}$^3$
	& 0.75 		
	& $0.505, 0.75, 0.9 $\\
\texttt{size\_LBFGS}$^3$			
	& 3 	    
	& $1, 2, 3, 4, 5, 10, 15, 20$ \\
\texttt{eps\_div\_zero\_grad}$^4$ 				
	& $10^{-10}$ 
	& $10^{-5}, 10^{-8\dagger}, 10^{-10}, 10^{-12\dagger}, 
	10^{-15}$ \\
\texttt{log\_zero\_safeguard}$^4$
	& $10^{-16}$ 
	& $10^{-4\dagger}, 10^{-8}, 10^{-12\dagger}, 10^{-16}, 
	10^{-24\dagger}, 10^{-32}$\\
\texttt{eps\_active\_set}$^4$
	& PDNR: $10^{-3}$ 
	& $10^{-1}, 10^{-3}, 10^{-5\dagger}, 10^{-8\dagger}$ \\
	& PQNR: $10^{-8}$
	& $10^{-1}, 10^{-3}, 10^{-5\dagger},10^{-8\dagger}$ \\
\bottomrule
\end{tabular}
}%
\\[1.5pt]
{\footnotesize $^{1-4}$Intel platform used for experiments; $^\dagger$values 
evaluated on Intel platform only}
\end{table*}

\subsection{Experiments}
\label{sec:methods:experiments}
An individual experiment is a job $j$ on platform $m$ 
solving a PDNR/PQNR row subproblem for dataset $d$ with SparTen solver 
$s$, parameter $p$, parameter value $v$, and random initialization $r$; all 
remaining software parameters are fixed at the default values listed in 
Table~\ref{tab:methods:params}. Certain experiments denoted with a dagger$^\dagger$ were 
run only on Intel hardware due to limited resources associated with the other 
architectures; this accounts for the larger number of experiments reported for these 
platforms. We conducted tests on these values to provide better resolution of the impact of 
the parameter where nearby values---i.e., on the bounds of the test range---contained 
uncertainty in the results. Furthermore, we split up the experiments across the Intel 
platforms by parameter, running the full set of experiments across all parameter values and 
all random initializations on a single platform. The superscripts denoted for each 
parameter in the table denote the Intel platform number specified in 
Table~\ref{tab:methods:hw}. Since we report only the number of function evaluations and 
outer iterations in our results, we expect that running our experiments in this way has 
produced valid results.

In all experiments, we fit a 5-component CP decomposition using 
a tolerance of $10^{-4}$ (i.e., the value of $\tau$ in Equation (20) 
in~\cite{HaPlKo15}, the violation of the Karush-Kuhn-Tucker (KKT) conditions, 
used as the stopping criterion for the methods we explore here). Computation of 
a CP decomposition using PDNR or PQNR in SparTen requires an initial guess of the model 
parameters---i.e., initial values for $\boldsymbol{\mathscr{M}}$ 
in~\eqref{eq:cpmodel}---drawn from a uniform distribution in the range $[0,1]$. 
As such, all experiments were replicated using 30 random 
initializations. We report results on the amount of computation required for 
convergence (i.e., the number of evaluations of the negative log likelihood 
objective function, $f(\boldsymbol{\mathscr{M}})$, defined in Equation~(4) 
of~\cite{HaPlKo15}) and the quality of the solution (i.e., the value of the 
negative log likelihood objective function). As each of our experiments 
consists of 30 replicates (i.e., 30 random initializations) across three CPU 
architectures, we report sample means and 95\% confidence intervals (as defined 
in~\cite{discrete-event-sim}) when presenting statistical trends in the results.
\clearpage

\section{Results}
\label{sec:results}
\begin{table}[!ht]
\centering
\caption{Experiments run on the different datasets and 
hardware 
platforms.}
\label{tab:experiments:overview:jobs}
\resizebox{\columnwidth}{!}{%
\begin{tabular}{lllrrrrrr}
\toprule
CPU & Solver & Data & Planned & Collected & Canceled & Converged & Max 
Iterations & Missing \\
\midrule
\multirow{6}{*}{ARM}&\multirow{3}{*}{PDNR}
&{\it chicago}	&1110&1110&4.8\%&82.2\%&13.0\%&0.0\%\\
& &{\it lbnl}	&1110&1110&10.5\%&76.5\%&13.0\%&0.0\%\\
& &{\it nell}	&1110&390&5.4\%&39.2\%&55.4\%&64.9\%\\
\cmidrule{2-9}
&\multirow{3}{*}{PQNR}
&{\it chicago}	&990&281&0.0\%&55.5\%&44.5\%&71.6\%	\\
& &{\it lbnl}	&990&237&0.0\%&0.0\%&100.0\%&76.1\%	\\
& &{\it nell}	&990&390&23.3\%&0.3\%&76.4\%&60.6\%	\\
\midrule
\multirow{6}{*}{IBM}
&\multirow{3}{*}{PDNR}
&{\it chicago}	&1110&855&5.4\%&77.8\%&16.8\%&23.0\%	\\
& &{\it lbnl}	&1110&692&11.3\%&73.3\%&15.4\%&37.7\%	\\
& &{\it nell}	&1110&424&51.2\%&12.0\%&36.8\%&61.8\%	\\
\cmidrule{2-9}
&\multirow{3}{*}{PQNR}
&{\it chicago}	&990&676&10.2\%&76.3\%&13.5\%&31.7\%	\\
& &{\it lbnl}	&990&293&61.8\%&0.0\%&38.2\%&70.4\%	\\
& &{\it nell}	&990&481&31.0\%&6.6\%&62.4\%&51.4\%	\\
\midrule
\multirow{6}{*}{Intel}
&\multirow{3}{*}{PDNR}
&{\it chicago}	&1680&1673&5.0\%&86.4\%&8.6\%&0.4\%	\\
& &{\it lbnl}	&1680&1663&11.0\%&80.6\%&8.4\%&1.0\%\\
& &{\it nell}	&1680&1643&44.7\%&42.2\%&13.1\%&2.2\%	\\
\cmidrule{2-9}
&\multirow{3}{*}{PQNR}
&{\it chicago}	&1440&1434&12.1\%&78.6\%&9.3\%&0.4\%	\\
& &{\it lbnl}	&1440&1363&78.0\%&0.0\%&22.0\%&5.3\%	\\
& &{\it nell}	&1440&1424&68.8\%&10.1\%&21.1\%&1.1\%	\\
\bottomrule
\end{tabular}%
}
\end{table}
In this section we analyze the results of the parameter 
sensitivity experiments and describe the statistical 
relationships between the convergence properties of the 
PDNR 
and PQNR methods and their input parameters.

In total, 21,960 unique experiments were planned, 
accounting 
for running PDNR and PQNR with random initializations 
across 
all parameter value ranges on the various hardware 
architectures described in 
Sections~\ref{sec:methods:hardware} 
and~\ref{sec:methods:params:description}. An experiment 
{\it 
converged} if the final KKT violation is less than the 
value of 
$\tau = \num{e-4}$; an experiment reached {\it maximum 
iterations} if the number of outer iterations exceeded 
the 
maximum limit (i.e., \texttt{max\_outer\_iterations}) 
and did 
not converge; an experiment was {\it canceled} if it 
exceeded 
the wall-clock limit (i.e., SparTen neither converged 
to a 
solution nor reached maximum number of outer iterations 
within 
12 hours); and an experiment was {\it missing} if it 
did not 
run due to a failure of the system to launch the 
experiment or 
other system issue. Of the planned experiments, we 
collected 
data from 16,139 experiments.

Table~\ref{tab:experiments:overview:jobs} presents the 
number of experiments planned as defined above and the number 
of planned experiments where data was collected (i.e., planned 
minus missing). For those collected, the table shows the 
percentage of experiments that were canceled, converged, or exceeded the 
maximum iterations. We note that the most complete set of experiment results 
were obtained on the Intel platforms. Although there are many missing 
experiment results for the IBM and ARM platforms, we attempt to identify 
patterns in the data we collected if there is strong evidence to support our 
claims. We note that a few parameters (\texttt{eps\_active\_set}, 
\texttt{min\_variable\_nonzero\_tolerance}, \texttt{suff\_decrease\_tolerance}, 
\texttt{damping\_increase\_tolerance}, \texttt{damping\_decrease\_tolerance}) 
showed no 
statistically significant differences across the range of input values used in 
the experiments. We conjecture that we did not find values where the parameters 
display sensitivities in the chosen tensor problems, thus it remains unclear if 
this behavior holds in general.

\subsection{General Convergence Results on Real-World Data}

\begin{figure}[ht!]
	\centering
	\includegraphics[width=\textwidth]{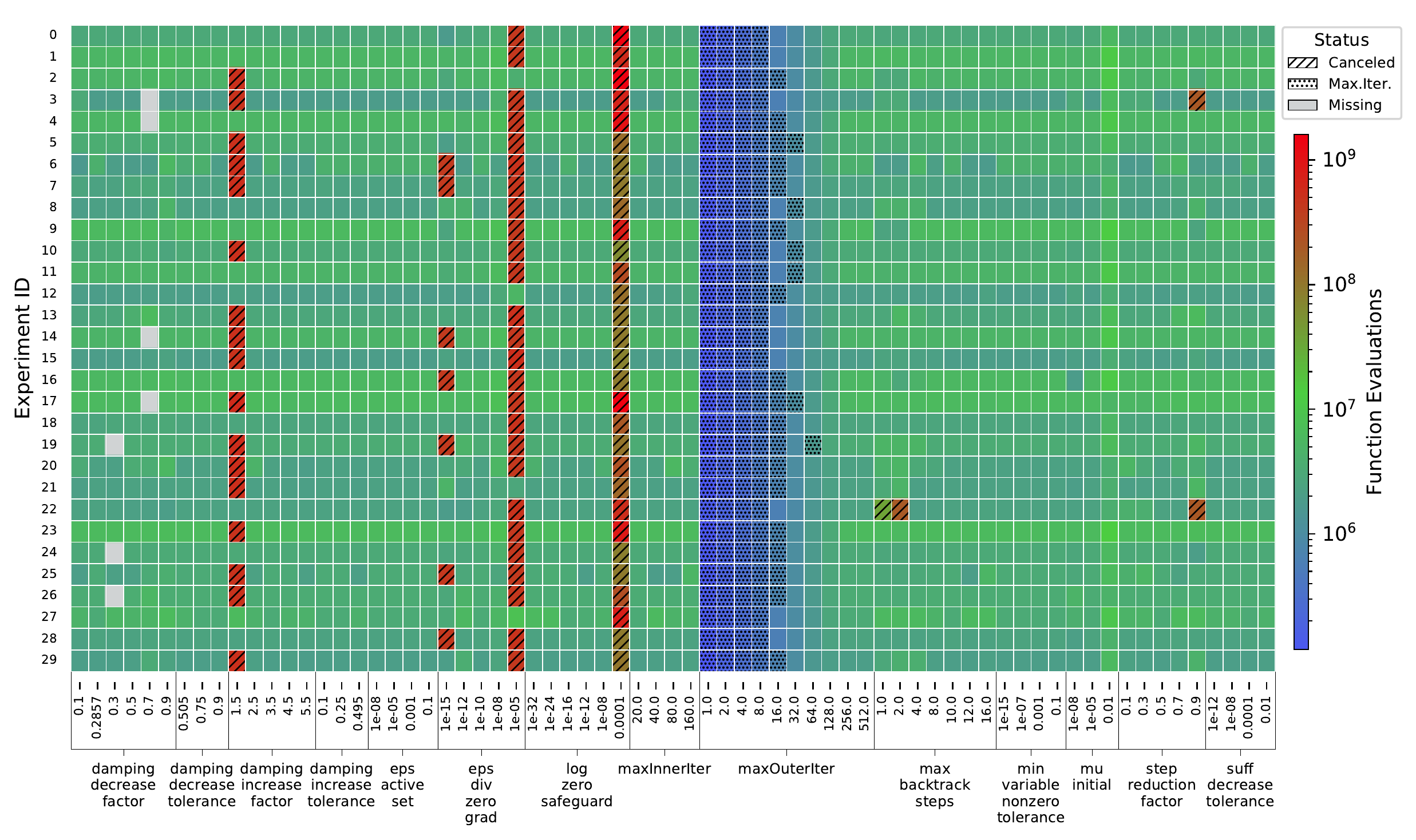}
	\caption{Example heatmap illustrating results.}
	\label{fig:results:heatmap:example}
\end{figure}

To illustrate general convergence behavior, we present the results of the 
experiments in a heatmap for a given dataset, method, and hardware platform.
Figure~\ref{fig:results:heatmap:example} presents an example heatmap, 
where each square represents the total number of objective function evaluations 
of an experiment, with random initializations across the rows and parameter 
values used (with all other parameter values set to their default values) 
across the columns. The complete set of heatmaps for all experiments can be 
found in Appendix~\ref{app:supplementary:heatmaps}. These results illustrate 
that there are certain ranges of parameter values that 
lead to good or bad convergence behaviors in general. 

The colours in the heatmaps denote the outcomes of the experiments as follows. 
Green shades are consistent with {\it converged} experiments. Vertical bands 
not shaded green identify values that may impact algorithm performance, due 
either to iteration constraints (blue hues) or excessive computations 
corresponding to slow convergence or stagnation (red hues). Hatches denote 
non-convergent exit status. Grey represents {\it missing} data, i.e., 
experiments that were planned but never conducted due to resource 
limitations---e.g., dequeued by the cluster administrator---or a system 
failure. Solid columns of a single shade indicating the same convergence 
behavior across all 30 random initializations. Nearly solid column lines of the 
same shade indicate similar behavior, but also that there is some sensitivity 
of those parameter values to the initial starting point of the iterative 
methods. 

\textbf{Observation 1: Convergence properties are demonstrated empirically.}
As discussed in Section~\ref{sec:intro}, applying PDNR 
and PQNR to real-world data has been explored previously in the literature only 
for a single problem. From Table~\ref{tab:experiments:overview:jobs}, we 
observed that PQNR is canceled more than PDNR in the allotted time across 
datasets and CPU platforms. This confirms our intuition, since it is a 
classical result in iterative methods that damped Newton methods converge 
quadratically, in comparison to quasi-Newton methods, which converge 
superlinearly. Specifically, PQNR calls the objective function $2.7$ times more 
than PDNR on average on the {\it chicago} data and fails to converge in the 
allotted time for any experiment on {\it lbnl} data across all hardware 
platforms. By contrast, PDNR converges in 86\% of {\it lbnl} experiments across 
platforms when only 32 outer iterations are allowed. 

\textbf{Observation 2: There is consistent convergence behavior across many 	
ranges of parameter values.} Looking across the range of values for each 
parameter, there are many cases where there is a distinct change in 
behavior---e.g., from converged to canceled when  
\texttt{log\_zero\_safeguard}  is greater than $10^{-8}$ in PDNR experiments on 
{\it chicago} as shown in Figure~\ref{fig:results:heatmap:intel:chicago}. This 
distinct, repeatable behavior can be used to guide the choices for both a 
general set of parameter defaults and tuning of parameters for specific data.

\textbf{Observation 3: PDNR and PQNR are not necessarily sensitive to the same 
parameters.} In some cases PQNR converges where PDNR is canceled using the 
same random initial starting point. Several interesting patterns from the {\it 
nell} results should also be noted. First, the solution space is highly 
sensitive to random initialization (seen as hatched red bands in 
Figure~\ref{fig:results:heatmap:intel:nell}). Thus statements quantifying 
average behavior are more uncertain. Second, where PQNR does converge on {\it 
nell}, it is for extreme values where convergence is unlikely in other 
cases---for example, using the following parameter values:
\begin{itemize}
\item \texttt{eps\_div\_zero\_grad} $\leq 10^{-12}$,
\item \texttt{log\_zero\_safeguard} $\geq 10^{-8}$,
\item \texttt{max\_backtrack\_steps} $\leq 4$.
\end{itemize}

In the next section, we present analysis of the results that expands on these 
observations.

\subsection{Sensitivity of Convergence and Solution Behavior}
We are interested in comparing the convergence results across parameter values 
and solvers. In the following sections, we will refer to 
Tables~\ref{tab:results:chicago-crime-comm:fevals},
~\ref{tab:results:lbnl-network:fevals}, and
~\ref{tab:results:nell-2:fevals}, which compare the computational costs between 
PDNR and PQNR by reporting the mean number of objective function evaluations of 
{\it converged} jobs from the three data sets for each parameter value on all 
computer hardware platforms with 95\% confidence intervals presented as 
percentages above and below the mean. We proceed by analyzing the results by 
the categories of parameters described in 
Section~\ref{sec:methods:params:description}: 
general CP-APR parameters, line search parameters, damped Newton step 
parameters, Quasi-Newton step parameters, and numerical stability parameters.

\subsubsection{CP-APR}
As noted above, damped Newton methods converge faster than Quasi-Newton 
methods 
(quadratic convergence versus super-linear, respectively). Thus, we are 
interested in determining the impact on convergence of constraining the 
maximum 
number of allowable iterations in the outer and inner solvers. Note that in 
experiments that do not measure the effect of \texttt{max\_outer\_iterations} 
and \texttt{max\_inner\_iterations} explicitly, these values were fixed at 
100,000 and 20, respectively.

In all cases, there is a minimum value where PDNR converges to a solution. 
However, PQNR times out for all \texttt{max\_outer\_iterations} test values on 
{\it lbnl} and {\it nell}. Where comparisons can be made ({\it 
chicago}), PQNR calls the objective function $2.7$ times 
more than PDNR on average, although this relative cost between solvers 
decreases as \texttt{max\_outer\_iterations} grows. 

The effect of increasing the maximum number of inner iterations, 
\texttt{max\_inner\_iterations}, is similar. The non-monotonic increase in 
function 
evaluations for the maximum number of inner iterations can be explained by the 
trade-off between outer and inner iterations depending on a particular data 
set. The exception to this trend is the result that PDNR calls the objective 
function $2.33$ times more than PQNR on {\it nell} (see 
Table~\ref{tab:results:nell-2:fevals}), although this is most likely due to 
so few converged PQNR experiments and may not be statistically significant.

In principle, there is a value of parameter \texttt{eps\_active\_set} that 
will 
have an effect on convergence. In practice, however, we did not find that 
value 
for any dataset. We note that this result differs from~\cite{HaPlKo15}, which 
found setting the parameter to $10^{-3}$ leads to faster PDNR convergence, 
compared to a value of $10^{-8}$. Therefore, since algorithm sensitivity to 
this parameter is data-dependent, it is unreasonable to generalize behavior.

\subsubsection{Line search}
\begin{figure*}[!ht]
\centering
\begin{minipage}[t]{0.45\textwidth}
\includegraphics[width=\textwidth]{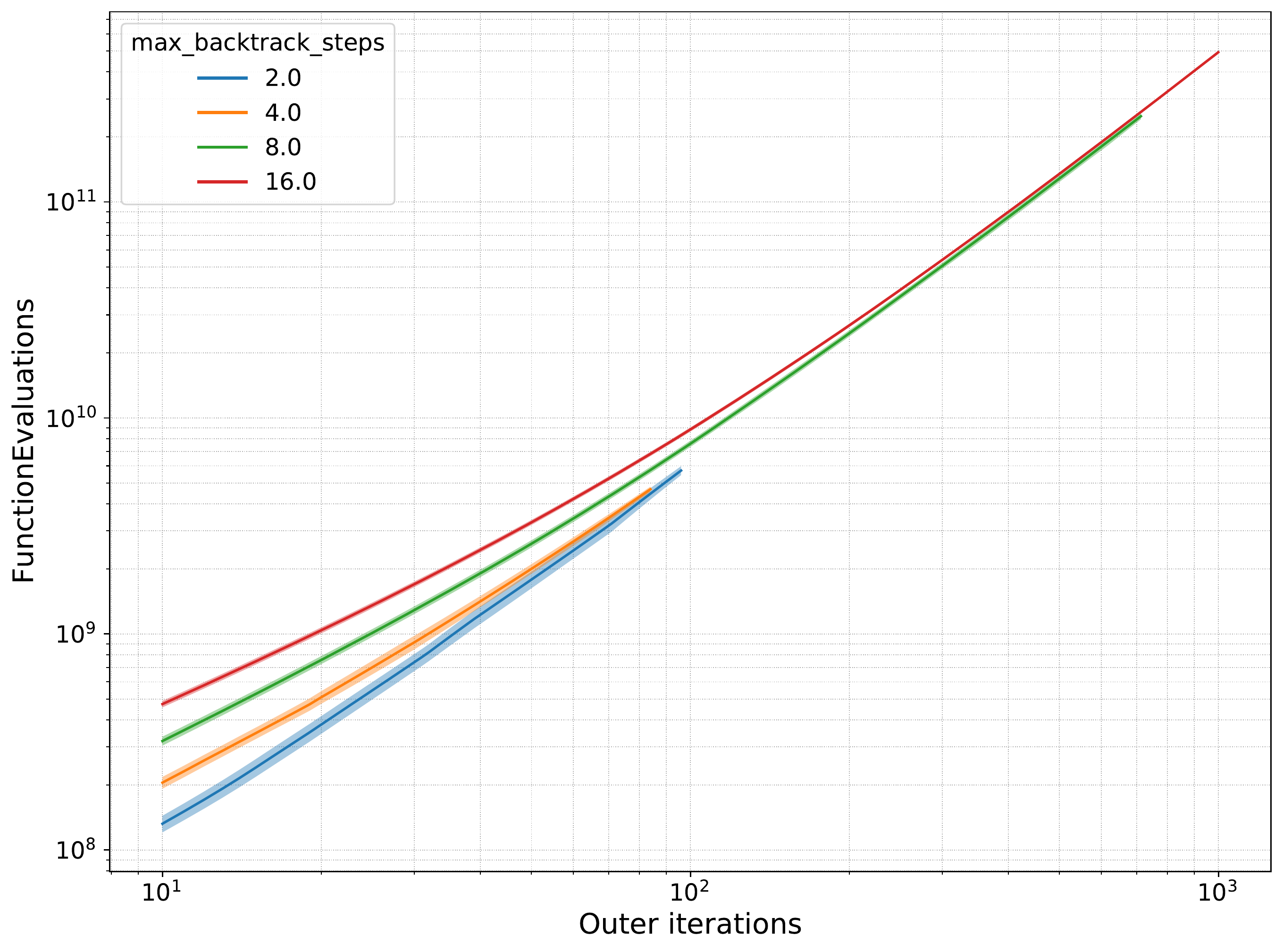}
\subcaption{The $x$-axis is truncated to emphasize the lower average cost when 
fewer \texttt{max\_backtrack\_steps} are allowed; the figure does not fully 
capture the high average cost when 16 maximum backtrack steps are allowed 
(PDNR, {\it lbnl}).}
\label{fig:results:uvs:pdnr:lbnl-network:max-backtrack-steps}
\end{minipage}
\quad
\begin{minipage}[t]{0.45\textwidth}
\includegraphics[width=\textwidth]{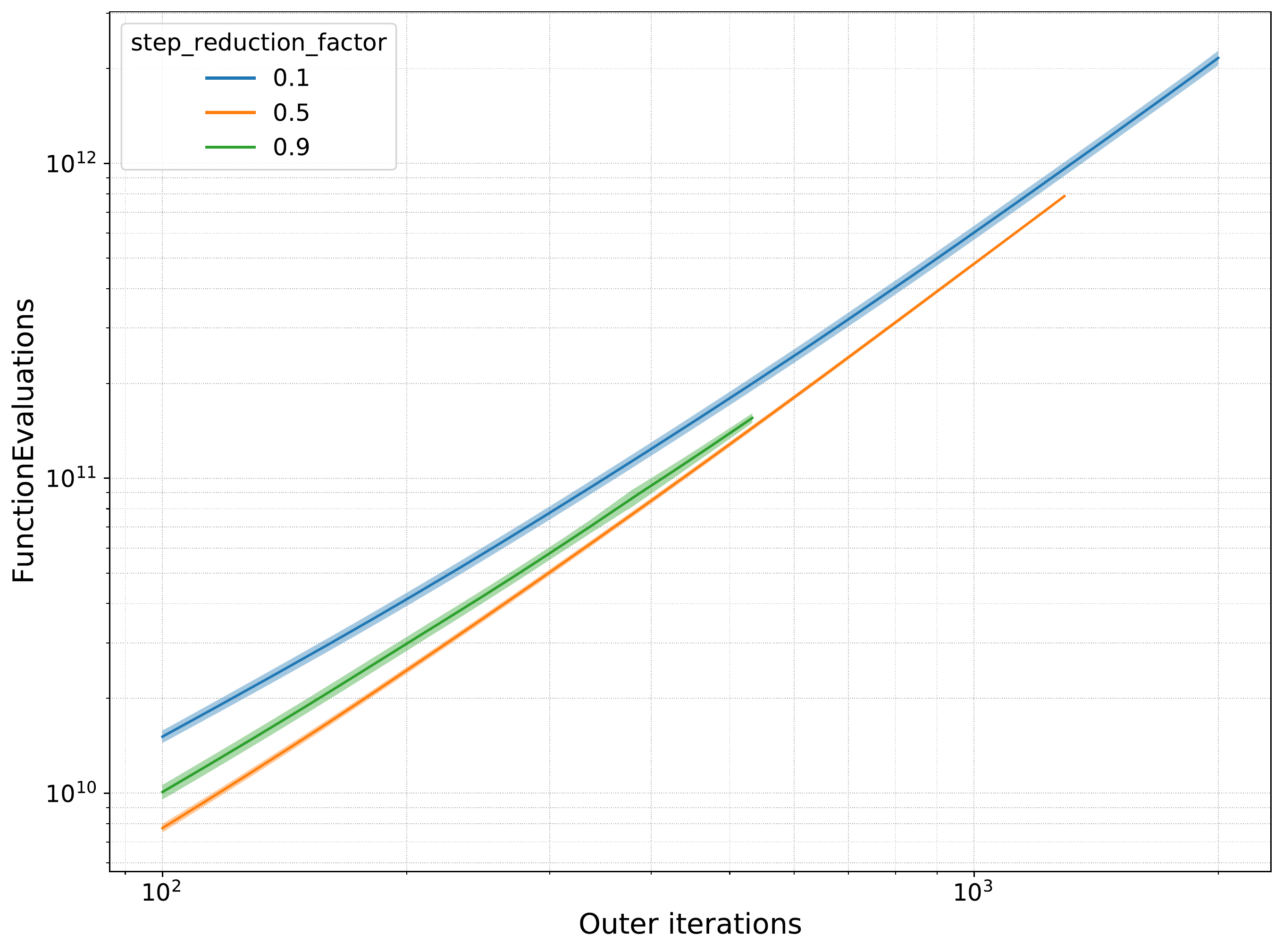}
\subcaption{The $x$-axis is truncated to demonstrate how high 	
\texttt{step\_reduction\_factor} may accelerate convergence on large, sparse 
tensor data (PDNR, {\it lbnl}).}
\label{fig:results:uvs:pdnr:lbnl-network:step-reduction-factor}
\end{minipage}
\end{figure*}

Allowing many backtracking steps during the line search, set by 
\texttt{max\_backtrack\_steps} may cause PDNR to waste effort; however, PQNR 
appears to perform better, in general, with more steps. PDNR is sensitive to 
the number of backtracking steps  on {\it chicago}: average work performed is 
less when the maximum number of allowed steps is large and more work is 
performed when the number of steps is small. On {\it lbnl}---the sparsest 
tensor problem considered---PDNR performs better with fewer backtracking steps 
(see Figure~\ref{fig:results:uvs:pdnr:lbnl-network:max-backtrack-steps}).

The line search parameter \texttt{step\_reduction\_factor} is used to 
reduce the line search step size between iterations ($\beta^{t+1}/\beta^t$ in 
Equation (17) in~\cite{HaPlKo15}).
On a large, sparse tensor problem, increasing this 
parameter may accelerate convergence. On the other hand, a small value makes 
convergence less certain. 
Figure~\ref{fig:results:uvs:pdnr:lbnl-network:step-reduction-factor}
illustrates this behavior on the {\it lbnl} data: the average 
total cost decreased by 77\% as \texttt{step\_reduction\_factor} increased from 
0.1 to 0.5 (SparTen default) and decreased another 28\% from 0.5 to 0.9. On 
{\it nell} data, PQNR {\it only} converged for large values (0.7, 0.9).

Threshold parameter \texttt{min\_variable\_nonzero\_tolerance} guarantees that 
the final Newton step length is nonzero. In theory, if this value is too 
large, the next iterate may overstep important information, forcing additional 
iterations to correct the misstep. On the other hand, if the value is too 
small, the algorithm may converge too slowly. We observed no statistically 
significant differences in PDNR algorithm performance when varying this 
parameter on {\it chicago}. The same appears true for PQNR, although 
we caution that empirical data is limited to results on {\it chicago} only. 
Addressing the result that no experiments converged for 
\texttt{min\_variable\_nonzero\_tolerance} \num{e-1} or \num{e-15} on 
{\it nell}, it is unlikely that size and sparsity play a role in parameter 
convergence; PDNR algorithm performance shows no significant difference 
on {\it lbnl}, the largest and sparsest tensor problem. 

The parameter \texttt{suff\_decrease\_tolerance} is used to assess if 
sufficient decrease in the objective function has been achieved in the line 
search---i.e., it is used to determine when to stop the line search 
iterations. 
A characterization virtually identical to that made for parameter 
\texttt{min\_variable\_nonzero\_tolerance} can also be stated here. In short, 
there is no significant performance difference when varying this parameter on 
PDNR or PQNR, although PQNR results are limited to {\it chicago} 
only.

\subsubsection{Damped Newton (PDNR)}
In this section, we discuss results of varying parameters that are used only 
by 
the PDNR solver. The damped 
Newton parameters control updates to the damping factor for the next iterate 
($\mu_k$ in Equation (16) in~\cite{HaPlKo15}). The damping parameter $\mu_k$ 
shifts the eigenvalues of the Hessian matrix, forcing it to be positive 
semidefinite and guaranteeing that a solution exists. For every outer 
iteration, the damping factor is initialized to \texttt{mu\_initial} and 
updated using the following parameters:
\begin{enumerate}
\item \texttt{damping\_decrease\_factor}
\item \texttt{damping\_decrease\_tolerance}
\item \texttt{damping\_increase\_factor}
\item \texttt{damping\_increase\_tolerance}
\end{enumerate}

\begin{figure*}
\centering
\begin{subfloat}[{\it chicago}]
{\includegraphics[width=0.3\textwidth]{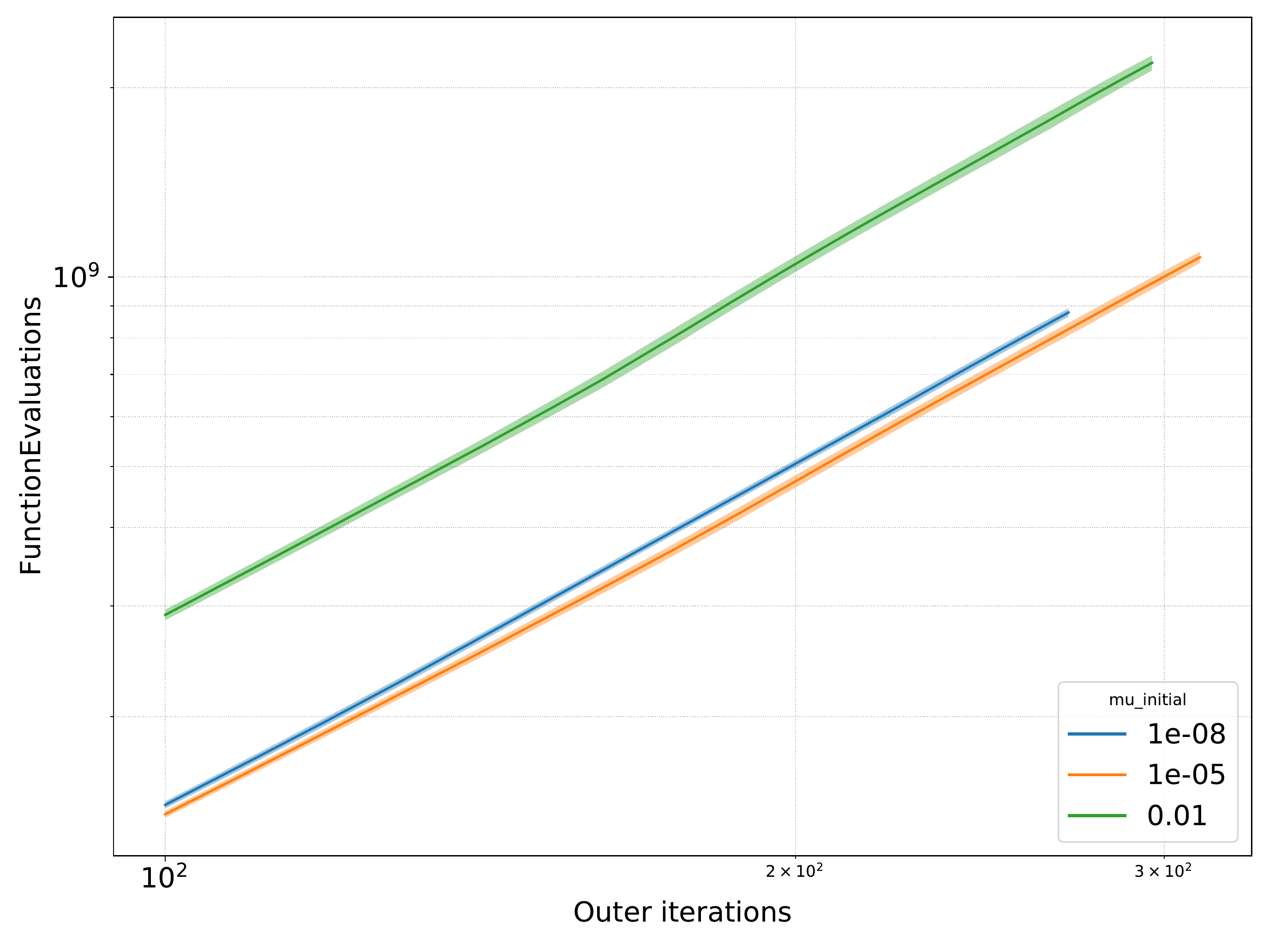}
\label{fig:results:uvs:pdnr:chicago:mu-initial}
}
\end{subfloat}%
\quad
\begin{subfloat}[{\it lbnl}]
{\includegraphics[width=0.3\textwidth]{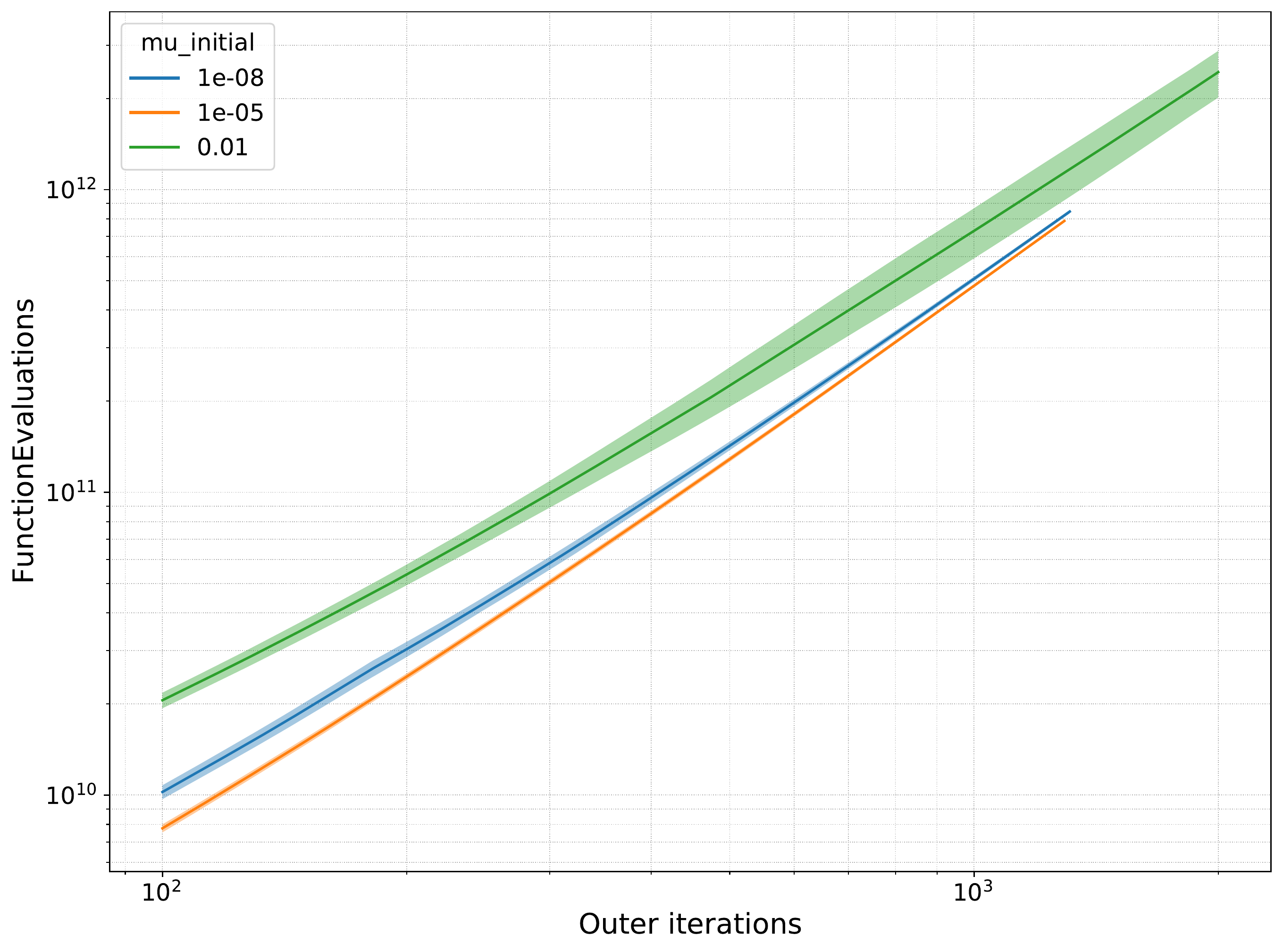}
\label{fig:results:uvs:pdnr:lbnl-network:mu-initial}
}
\end{subfloat}%
\quad
\begin{subfloat}[{\it nell}]
{\includegraphics[width=0.3\textwidth]{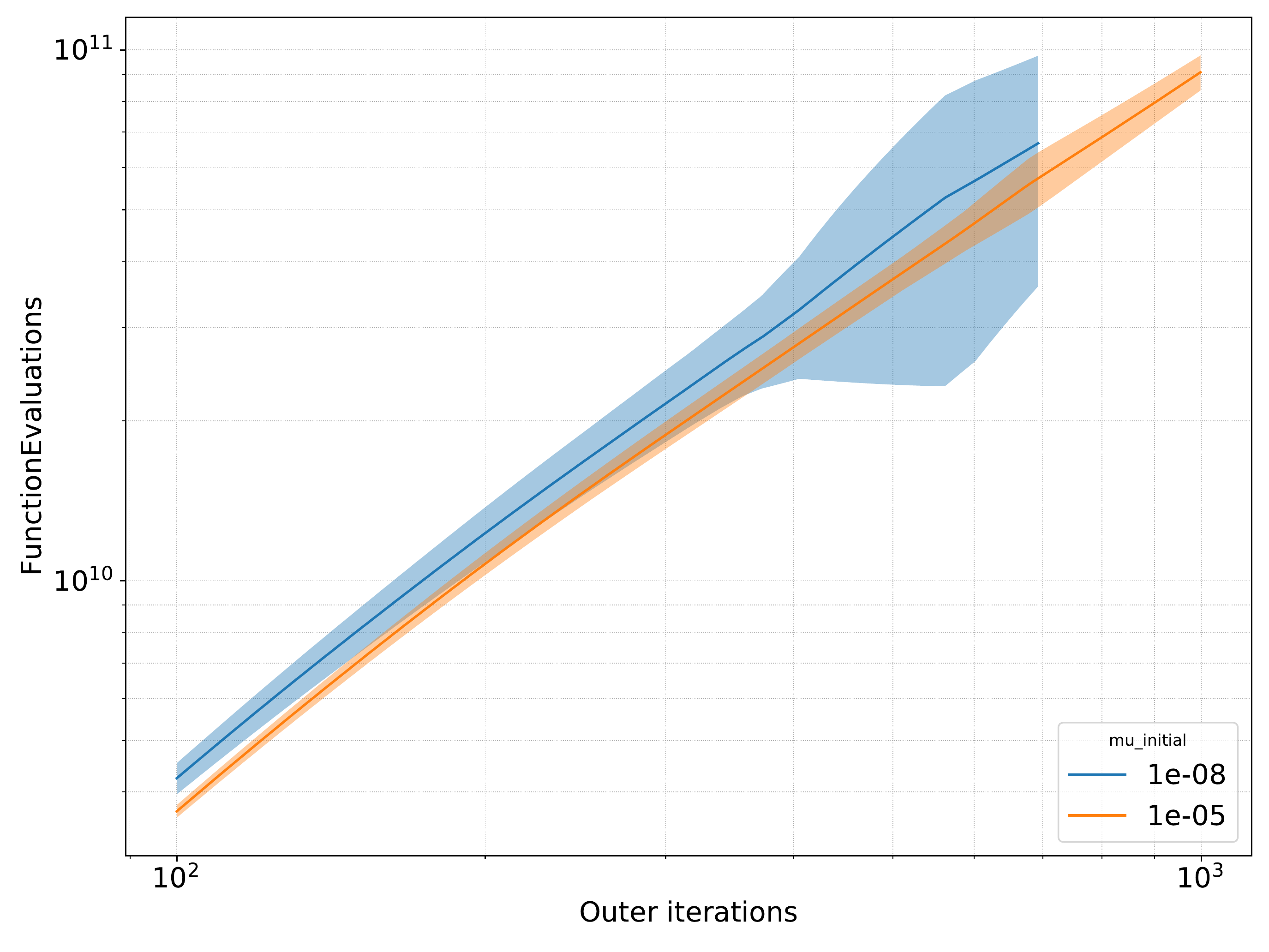}
\label{fig:results:uvs:pdnr:nell-2:mu-initial} 
}%
\end{subfloat}
\caption{PDNR, \texttt{mu\_initial}. Mean function evaluations and 
95\% CI. PDNR damps out Hessian information and is more prone to time-outs 
when \texttt{mu\_initial} is large (PDNR, {\it lbnl}).}
\label{fig:results:uvs:pdnr:mu-initial}
\end{figure*}

Hansen {\it et al.} predict in~\cite{HaPlKo15} that when the damping parameter 
$\mu$ is set too large, a loss of Hessian information follows, which impacts 
convergence:
\begin{quote}
``We expect larger values of $\mu_k$ to improve robustness by effectively 
shortening the step length and hopefully avoiding the mistake of setting 
too many variables to zero. However, a serious drawback to increasing 
$\mu_k$ is that it damps out Hessian information, which can hinder the 
convergence rate.''
\end{quote}
For example, when \texttt{mu\_initial} is large, the computational 
cost grows dramatically and time-outs become more likely, since the initial 
step length will at first be very small in every outer iteration and useful 
Hessian information is discarded in early stages of the inner loop solves. See 
Figure~\ref{fig:results:uvs:pdnr:mu-initial}.
Convergence is most likely for a large, but not too-large, value, i.e., 
\texttt{mu\_initial} $=10^{-5}$. Cost grows 177.2\% on {\it lbnl} and nearly 
doubles ($+92.2\%$) on {\it chicago} as \texttt{mu\_initial} grows from 
$10^{-5}$ to $10^{-2}$. It is important to note in the former case that this 
cost is skewed by one experiment that converged after nearly 42,000 outer 
iterations, in comparison to 1,300 for the other parameter values on average, 
illustrated in Figure~\ref{fig:results:uvs:pdnr:mu-initial}, where the $x$-axis 
is truncated to highlight the differences in total cost. Smaller values (i.e, 
$10^{-8}$) seem to perform better for {\it chicago}, the smallest, densest 
problem and larger values (i.e., $10^{-5}$) tend to perform better for large, 
sparse problems.

PDNR parameters \texttt{damping\_increase\_factor} and 
\texttt{damping\_decrease\_factor}, which control updates to the 
Hessian matrix damping parameter $\mu$, are two examples of algorithm 
parameters where convergence behavior is similar for values set within 
sensitivity constraints. SparTen rarely converges when the former is set too 
low (1.5); the likely effect is that the updated damping factor is insufficient 
to guarantee a well-conditioned Hessian and too many unimportant directions are 
considered when computing the search direction. Above the 1.5 bound, the cost 
in objective function calls does not change significantly. 

In our experiments, we found that varying either 
\texttt{damping\_decrease\_tolerance} or \texttt{damping\_increase\_tolerance} 
has no noticeable effect on either the number of calls to the objective 
function nor the number of outer iterations performed.

\subsubsection{Quasi-Newton (PQNR)}
In this section, we discuss varying parameter \texttt{size\_LBFGS}, which is 
the only software parameter used by the PQNR solver. PQNR uses a limited memory 
BFGS (L-BFGS) approach to approximate the inverse Hessian matrix in the 
Quasi-Newton step, with $M$ update pairs stored. The value of $M$ is set by the 
parameter \texttt{size\_LBFGS}. See~\cite{HaPlKo15} for algorithm details. More 
update pairs $M$ should provide a higher resolution to approximate the 
inverse Hessian. Intuitively, too few update pairs seems insufficient to 
compute an acceptable approximation. The only observable difference occurs when 
the update size is 1, using only the current iterate in the BFGS update.

\subsubsection{Numerical stability}
The numerical stability parameters \texttt{eps\_div\_zero\_grad} and 
\texttt{log\_zero\_safeguard} described in this Section are offset tolerances 
to avoid divide-by-zero floating point errors tailored to SparTen's C++ 
implementation. They do not appear in the corresponding Matlab Tensor Toolbox 
method \verb|cp_apr|. Their impact on convergence was consistent across 
combinations of solver, data, and CPU hardware.

The parameter \texttt{eps\_div\_zero\_grad} is an offset to avoid 
divide-by-zero floating point errors when computing the gradient and Hessian 
(Equation (10) in~\cite{HaPlKo15}). When \texttt{eps\_div\_zero\_grad} is 
large, gradient directions that do not lead to objective function improvements 
may be scaled the same as gradient directions that do lead to such 
improvements. Furthermore, the corresponding eigenvalues of the Hessian matrix 
are amplified and Hessian information may be lost when determining the next 
iterate. For example, PDNR loses Hessian information as 
\texttt{eps\_div\_zero\_grad} increases on {\it chicago} data; PDNR rarely 
converges and PQNR never converges when this parameter is relatively 
large---i.e. $10^{-5}$. Moreover, both algorithms are sensitive to the 
parameter's lower bound, as small values may be insufficient to avoid an 
ill-conditioned Hessian matrix. In either case, additional iterations follow to 
correct errors incurred by \texttt{eps\_div\_zero\_grad} values, large and 
small.

\begin{figure}[!ht]
\centering\includegraphics[width=.75\textwidth]{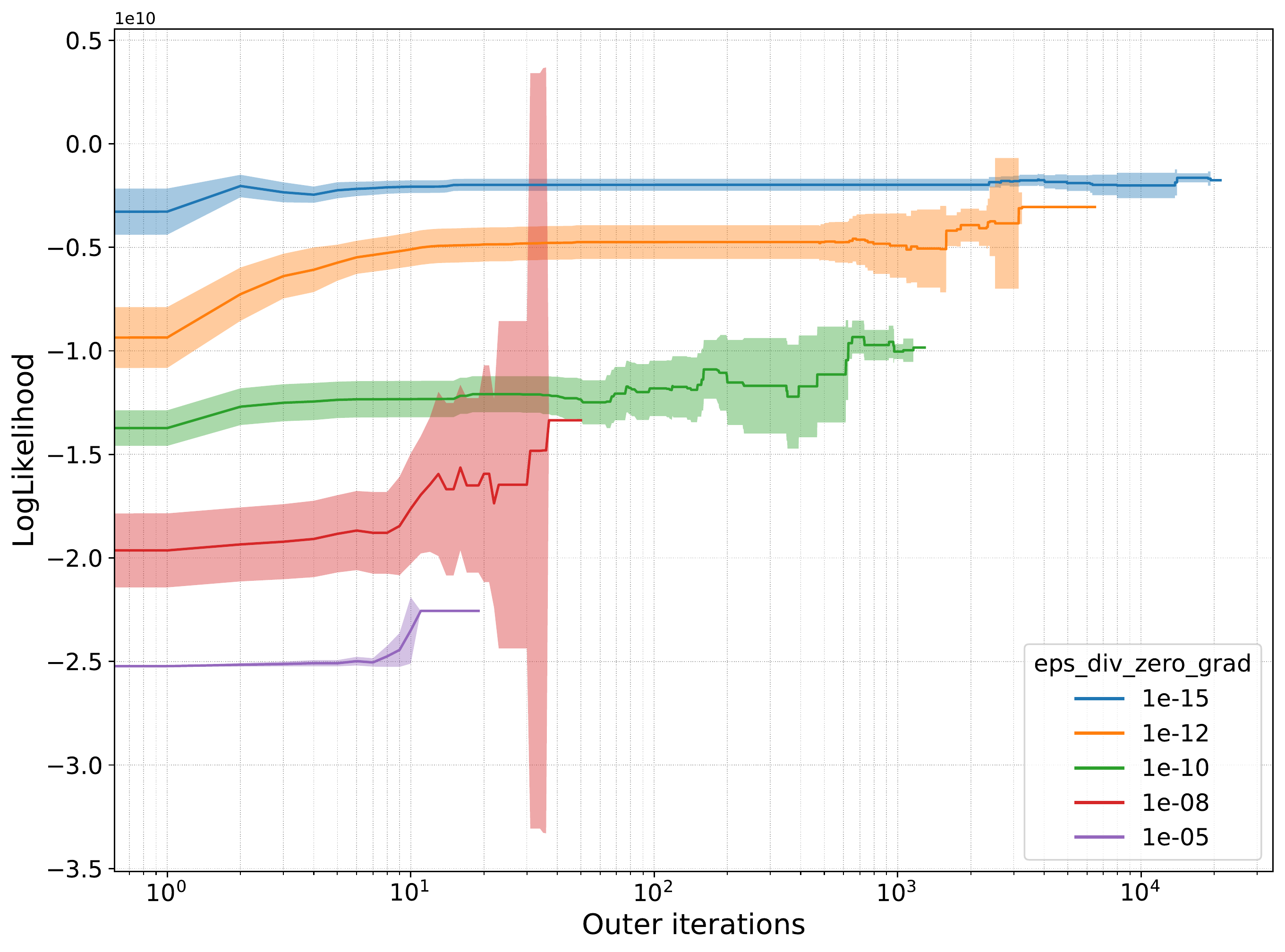}
\caption{Mean objective function values with 95\% confidence interval, varying 
\texttt{eps\_div\_zero\_grad} (PDNR, {\it lbnl}).}
\label{fig:results:uvs:pdnr:lbnl-network:eps-div-zero-grad:obj}
\end{figure}

Parameter sensitivities affect not only convergence 
behavior, but may also produce qualitatively different results. 
Figure~\ref{fig:results:uvs:pdnr:lbnl-network:eps-div-zero-grad:obj}
illustrates the effect where large \texttt{eps\_div\_zero\_grad}---and 
consequently, small step length---minimizes calls to the objective function 
{\it and} results in minimal objective function value. Most striking is 
that larger \texttt{eps\_div\_zero\_grad} decreases the objective function more 
than an order of magnitude. This result was collected from 79 of 90 planned 
PDNR experiments on {\it lbnl}, and thus we consider this interesting 
effect worthy of further investigation.

PDNR typically does not converge for large \texttt{log\_zero\_safeguard} values 
on large tensor problems. This parameter sets a nonzero offset in logarithm 
calculations to avoid explicitly computing $\log(0)$. High precision in 
logarithm computations tends to ensure the objective function is minimized 
accurately. When the value is too large, the calculated logarithm may be too 
small, and more backtracking steps are required to sufficiently decrease the 
objective function in the line search routine, making time-outs more likely. On 
the other hand, the effect of the parameter on convergence is indistinguishable 
for values smaller than $10^{-8}$ across all experiments.
\clearpage

\section{Conclusions}
Using results from more than 16,000 numerical experiments on several hardware 
platforms, we presented experimental results that expand our understanding of 
average PDNR and PQNR convergence on real-world tensor problems. We have shown that 
when using PQNR to compute large tensor decompositions convergence is 
less-likely under reasonable resource constraints. We have shown that some 
software parameters are sensitive to bounds on values. Further, we showed 
that varying several parameters can dramatically impact algorithm performance, and in some 
cases, may produce qualitatively different results.

Future work may address the issue of stagnation in Newton optimization methods 
for CP decompositions. We showed examples where the solver converged to 
a solution slowly but within the allotted time of 12 hours. For those 
experiments that timed out, it is unknown whether SparTen would eventually 
converge to a solution or stagnate without making progress. We anticipate that 
stagnation could be determined if the objective function values converge to a 
statistical steady state without satisfying the convergence criterion. Future 
development of SparTen may include dynamic updates to algorithm parameters 
based on local convergence information. Lastly, future experiments could 
explore coupled sensitivities among algorithm parameters, as this work was 
limited to single parameter, univariate analyses. Understanding the nature of 
bivariate (or even more complex) relationships among parameters may better 
inform end-users when searching for optimal parameter choices to run the 
SparTen methods.
\label{sec:conc}
\clearpage

\clearpage
\addcontentsline{toc}{section}{References}
\bibliography{arXiv_ref}

\clearpage
\appendix
\section{Detailed Experiment Results}
This section provides additional results for the experiments described in this 
analysis.

\subsection{Heatmaps}\label{app:supplementary:heatmaps}
Below are the outcomes of the planned experiments described in this report presented as heatmap images. They are organized by hardware platform, data set and solver used in each experiments, as described in Section~\ref{sec:methods}.

\begin{figure}[!htp]
	\centering
	\caption{Experiment outcomes for {\it chicago-crime-comm} data on Intel platform.}
	\label{fig:results:heatmap:intel:chicago}
	\begin{subfigure}{\textwidth}
		\includegraphics[width=\textwidth]{{fig.heatmap.evals.intel.chicago-crime-comm.Damped-Newton}.pdf}
		\subcaption{{\it PDNR}}
		\label{fig:results:heatmap:intel:pdnr:chicago}
	\end{subfigure}
	~
	\begin{subfigure}{\textwidth}
		\includegraphics[width=\textwidth]{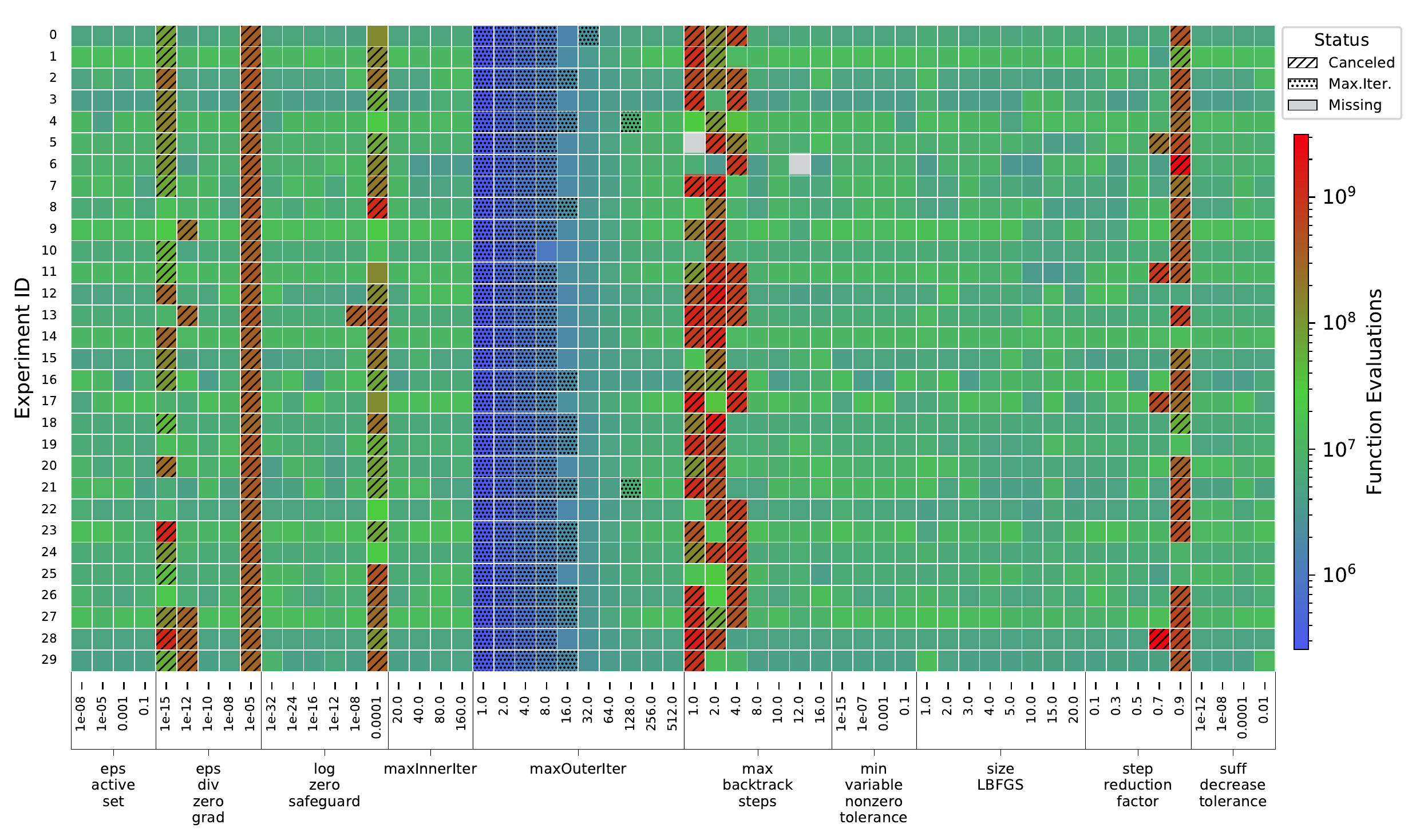}
		\subcaption{{\it PQNR}}
		\label{fig:results:heatmap:intel:pqnr:chicago}
	\end{subfigure}
\end{figure}
\clearpage
\begin{figure}[!htp]
	\centering
	\caption{Experiment outcomes for {\it lbnl-network} data on Intel platform.}
	\label{fig:results:heatmap:intel:lbnl}
	\begin{subfigure}{\textwidth}
		\includegraphics[width=\textwidth]{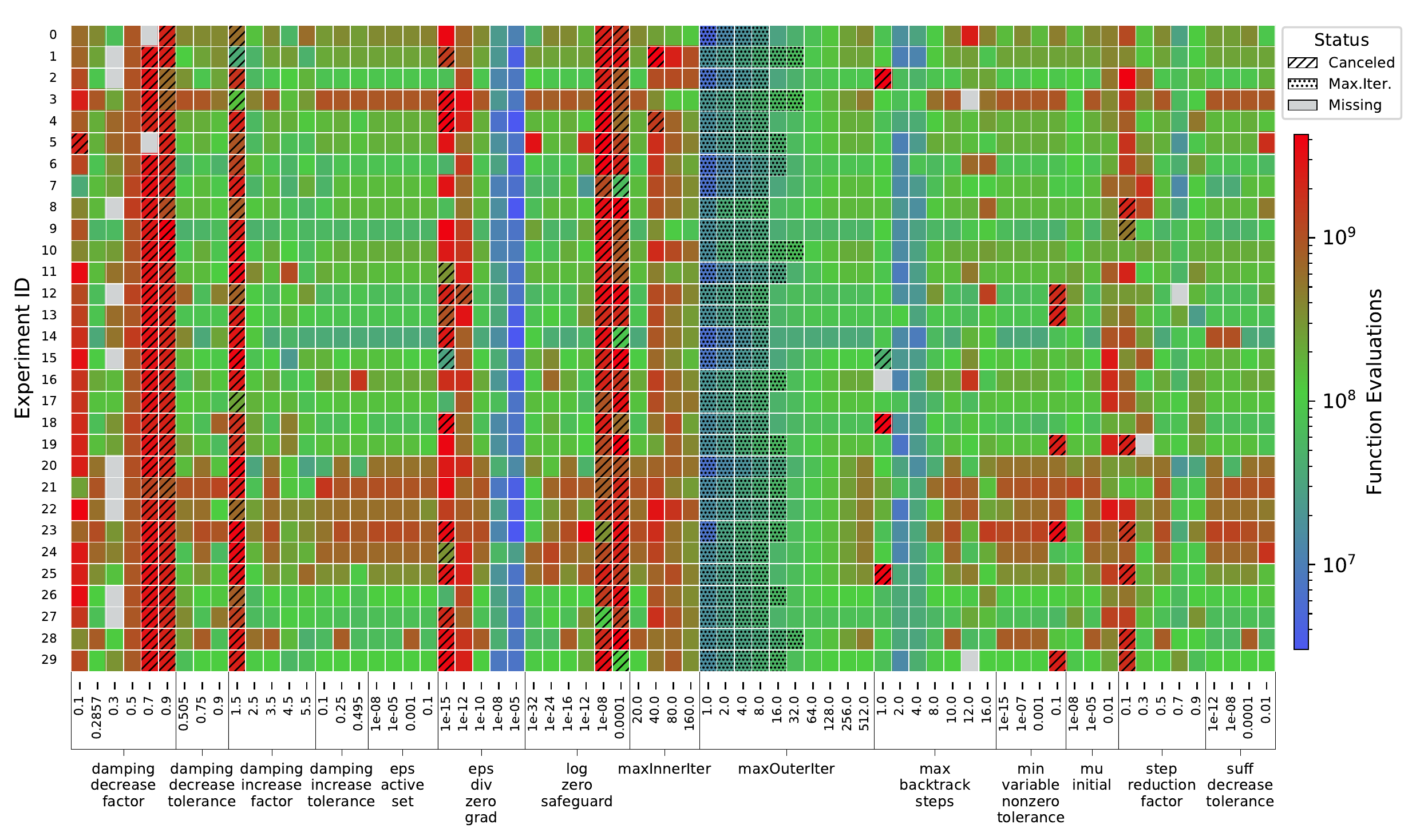}
		\subcaption{{\it PDNR}}
		\label{fig:results:heatmap:intel:pdnr:lbnl}
	\end{subfigure}
	~
	\begin{subfigure}{\textwidth}
		\includegraphics[width=\textwidth]{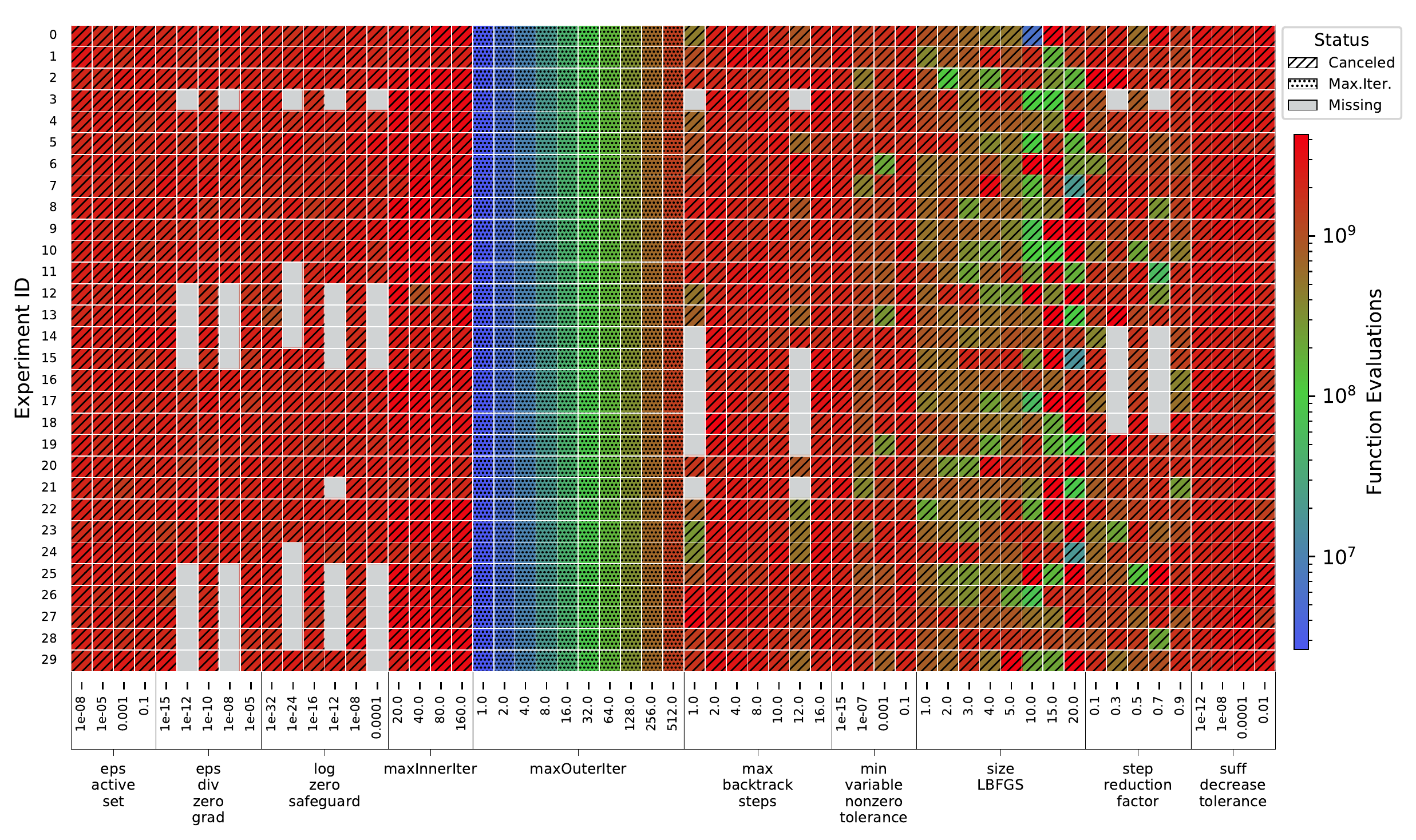}
		\subcaption{{\it PQNR}}
		\label{fig:results:heatmap:intel:pqnr:lbnl}
	\end{subfigure}
\end{figure}
\clearpage
\begin{figure}[!htp]
	\centering
	\caption{Experiment outcomes for {\it nell-2} data on Intel platform.}
	\label{fig:results:heatmap:intel:nell}
	\begin{subfigure}{\textwidth}
		\includegraphics[width=\textwidth]{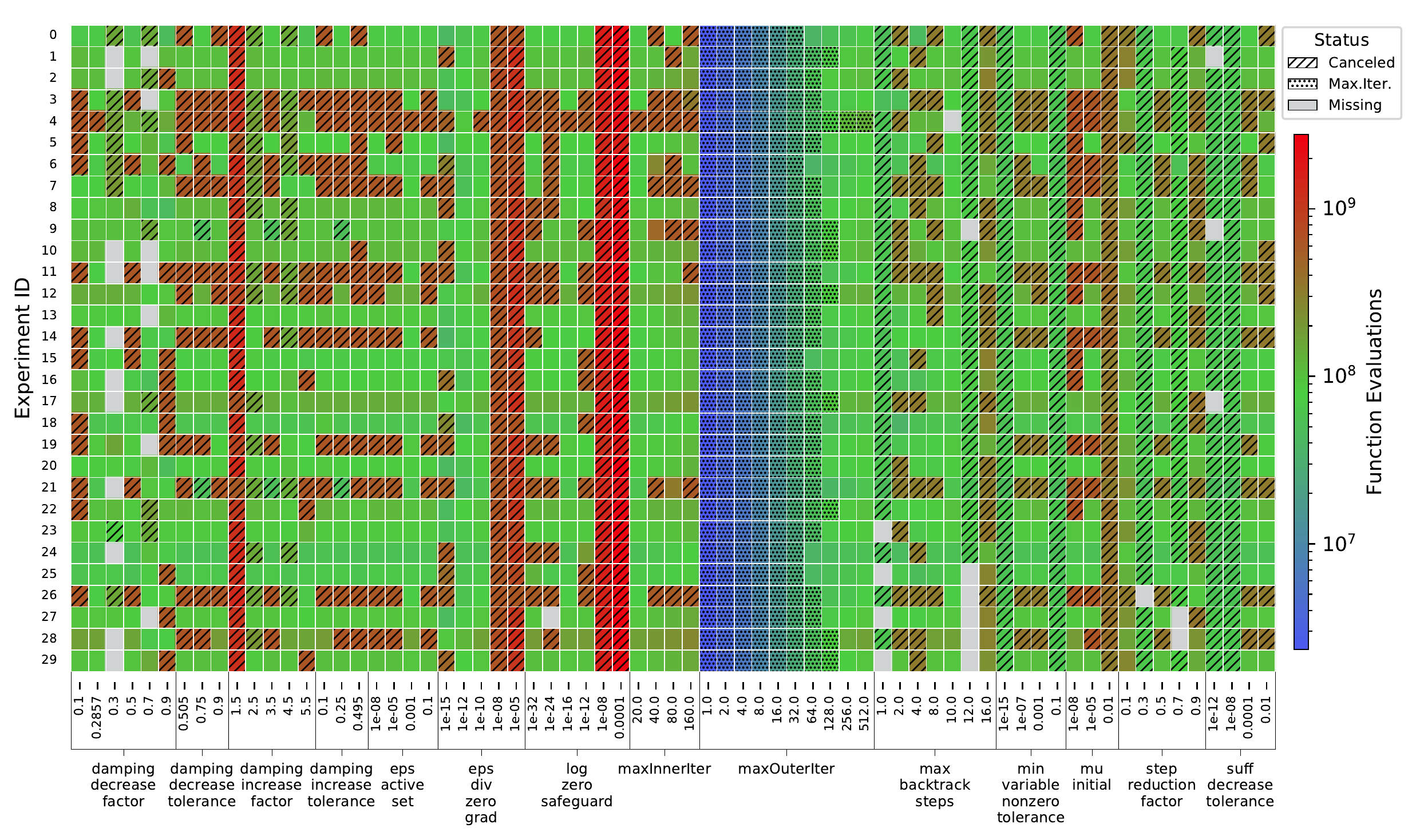}
		\subcaption{{\it PDNR}}
		\label{fig:results:heatmap:intel:pdnr:nell}
	\end{subfigure}
	~
	\begin{subfigure}{\textwidth}
		\includegraphics[width=\textwidth]{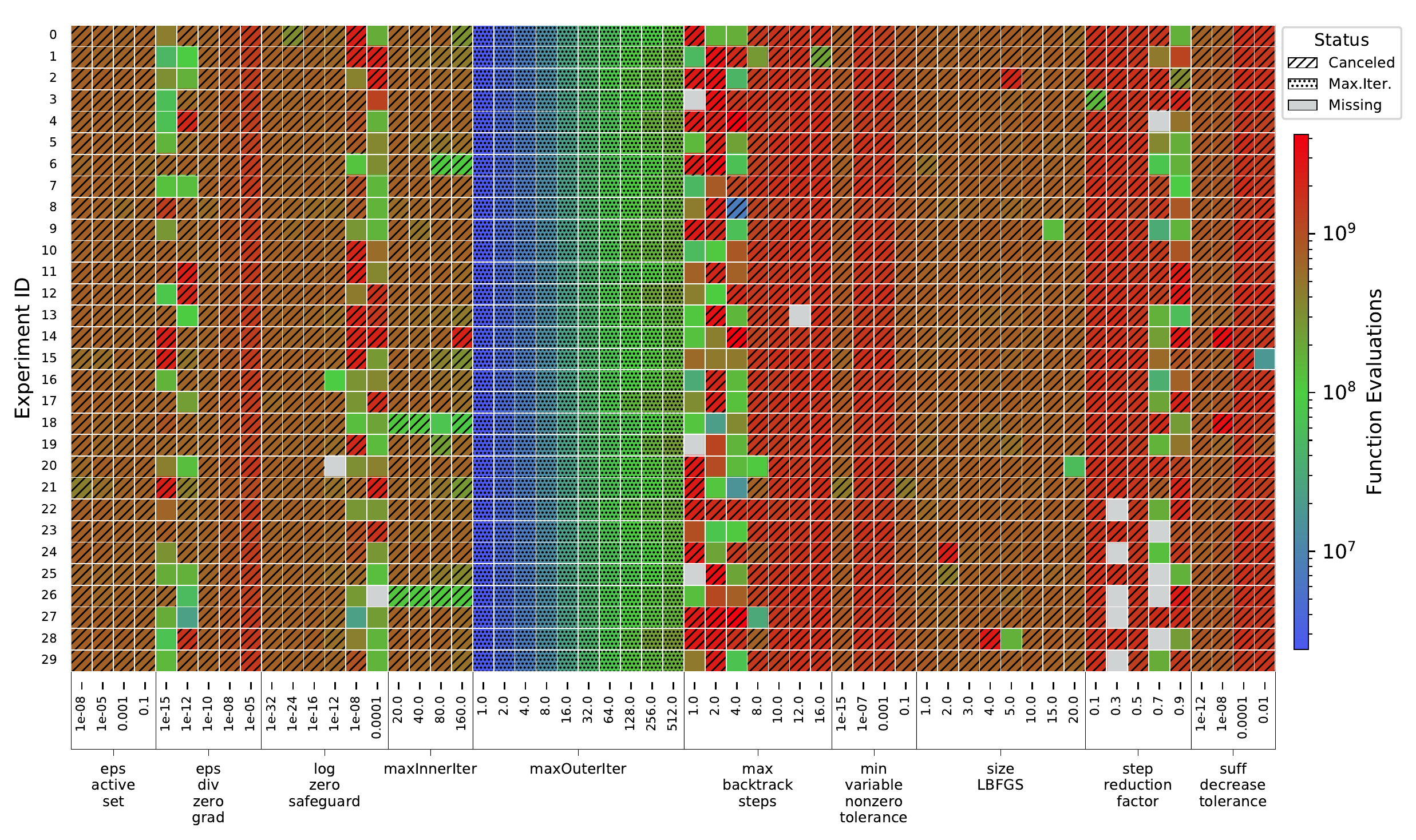}
		\subcaption{{\it PQNR}, {\it nell-2}}
		\label{fig:results:heatmap:intel:pqnr:nell}
	\end{subfigure}
\end{figure}
\clearpage

\begin{figure}[!htp]
	\centering
	\caption{Experiment outcomes for {\it chicago-crime-comm} data on ARM platform.}
	\label{fig:results:heatmap:arm:chicago}
	\begin{subfigure}{\textwidth}
		\includegraphics[width=\textwidth]{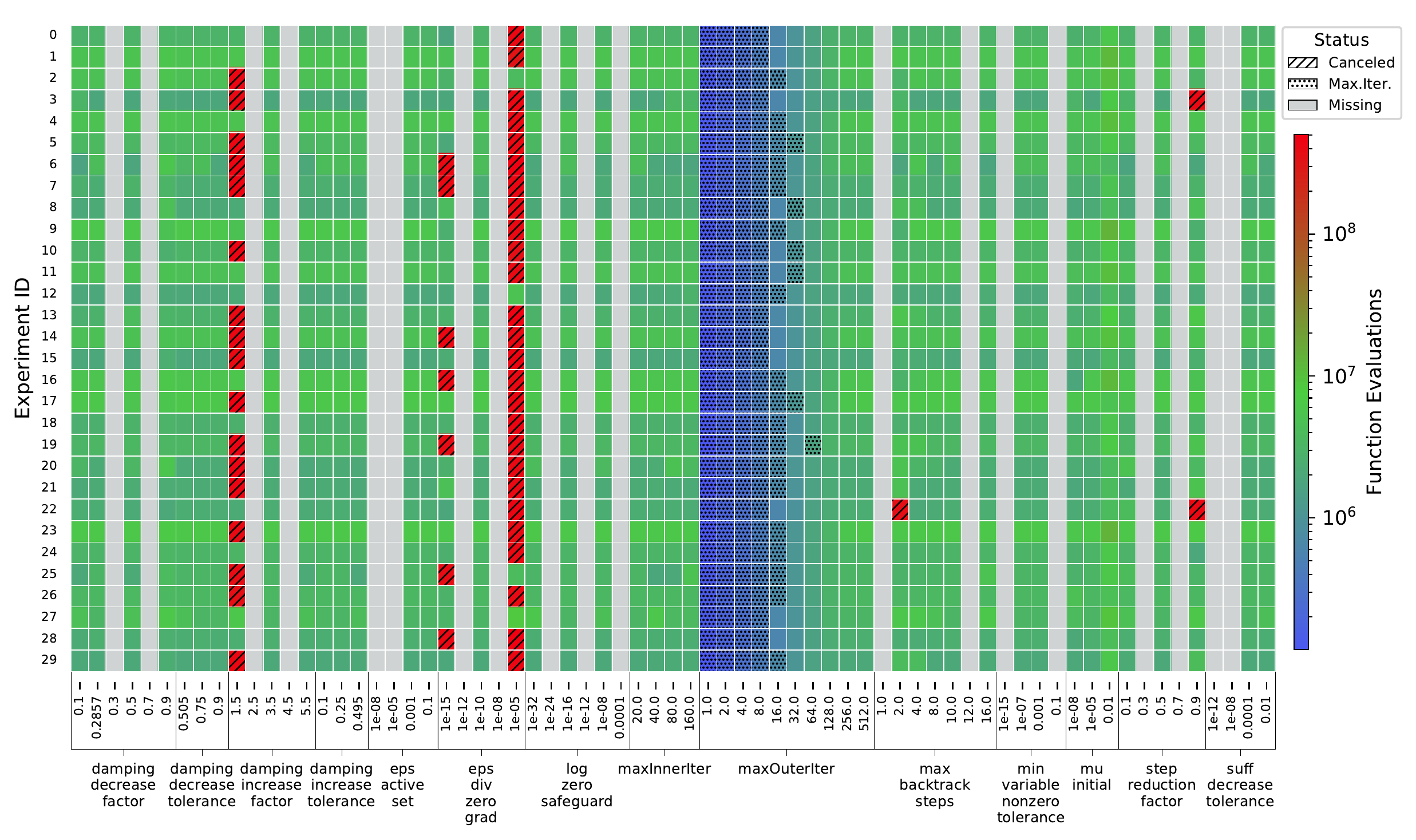}
		\subcaption{{\it PDNR}, {\it chicago-crime-comm}}
		\label{fig:results:heatmap:arm:pdnr:chicago}
	\end{subfigure}
	~
	\begin{subfigure}{\textwidth}
		\includegraphics[width=\textwidth]{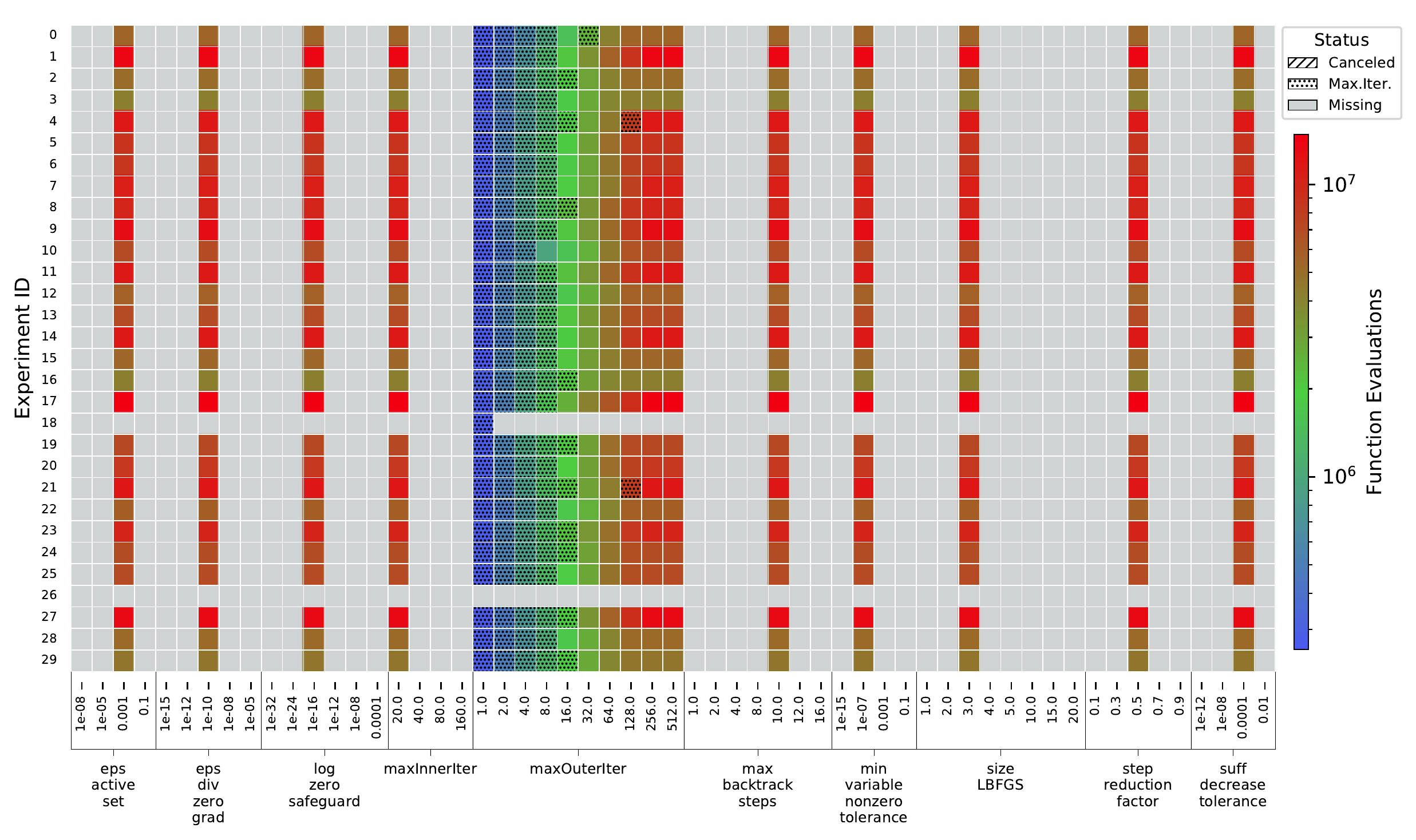}
		\subcaption{{\it PQNR}, {\it chicago-crime-comm}}
		\label{fig:results:heatmap:arm:pqnr:chicago}
	\end{subfigure}
\end{figure}
\clearpage
\begin{figure}[!htp]
	\centering
	\caption{Experiment outcomes for {\it lbnl-network} data on ARM platform.}
	\label{fig:results:heatmap:arm:lbnl}
	\begin{subfigure}{\textwidth}
		\includegraphics[width=\textwidth]{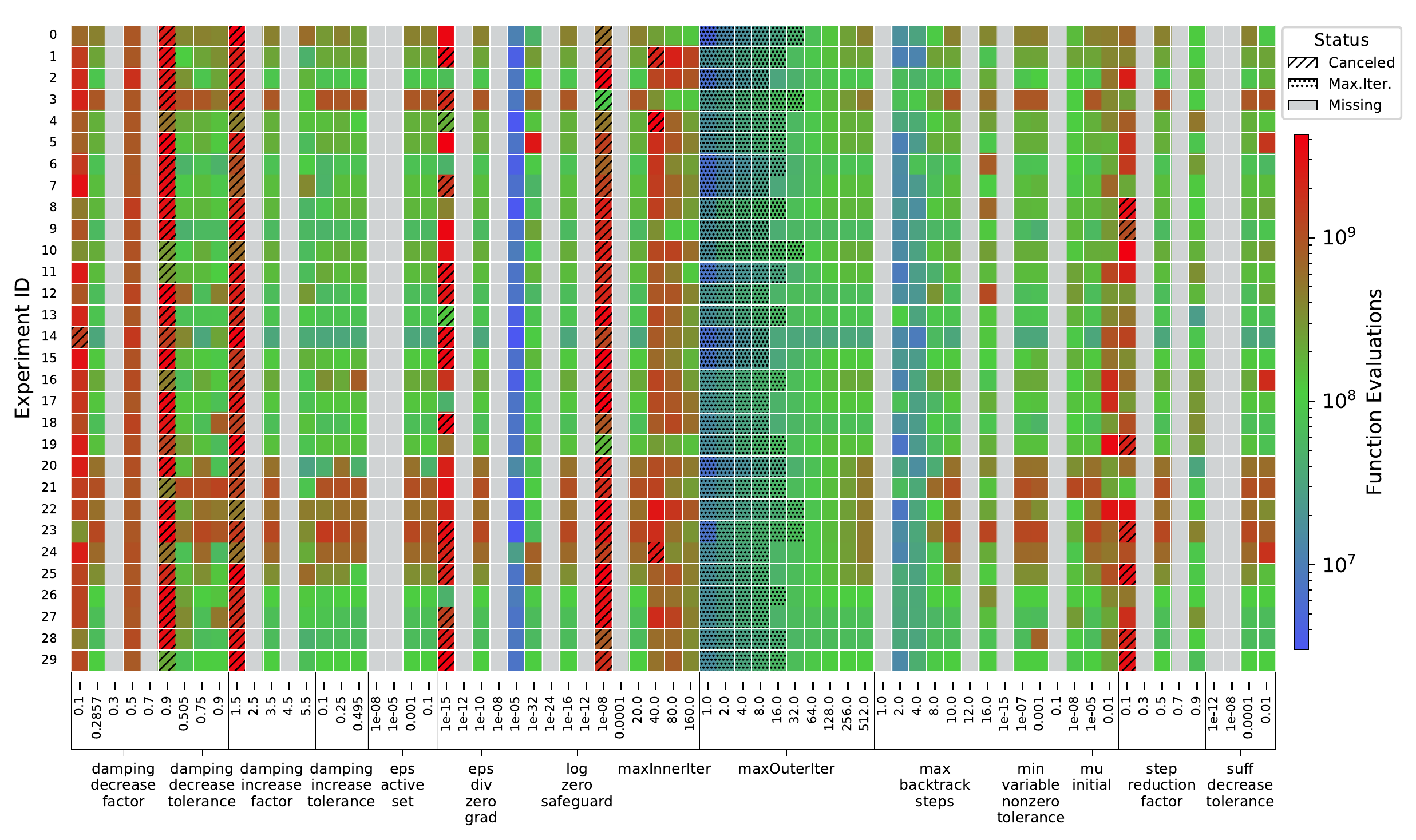}
		\subcaption{{\it PDNR}, {\it lbnl-network}}
		\label{fig:results:heatmap:arm:pdnr:lbnl}
	\end{subfigure}
	~
	\begin{subfigure}{\textwidth}
		\includegraphics[width=\textwidth]{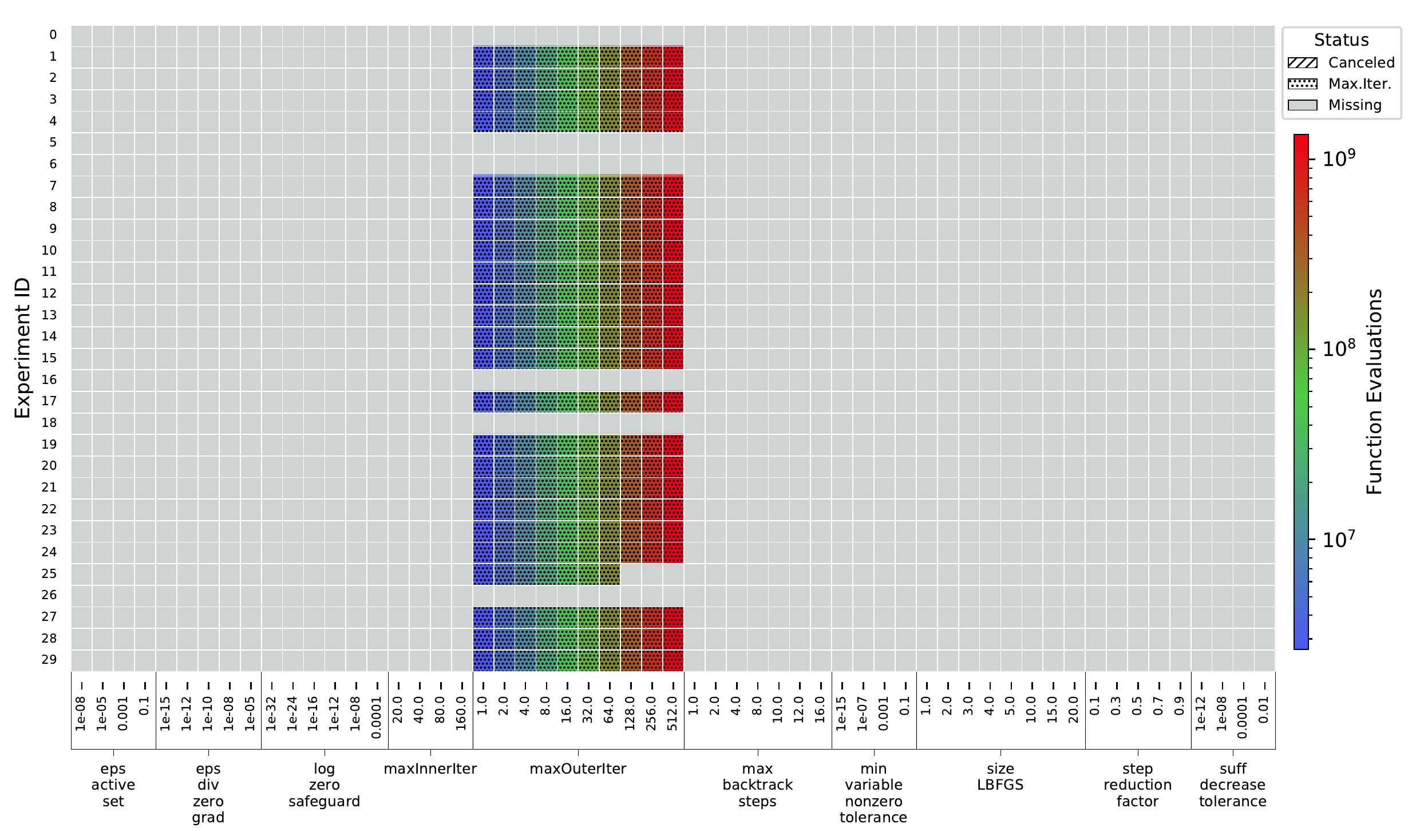}
		\subcaption{{\it PQNR}, {\it lbnl-network}}
		\label{fig:results:heatmap:arm:pqnr:lbnl}
	\end{subfigure}
\end{figure}
\clearpage
\begin{figure}[!htp]
	\centering
	\caption{Experiment outcomes for {\it nell-2} data on ARM platform.}
	\label{fig:results:heatmap:arm:nell}
	\begin{subfigure}{\textwidth}
		\includegraphics[width=\textwidth]{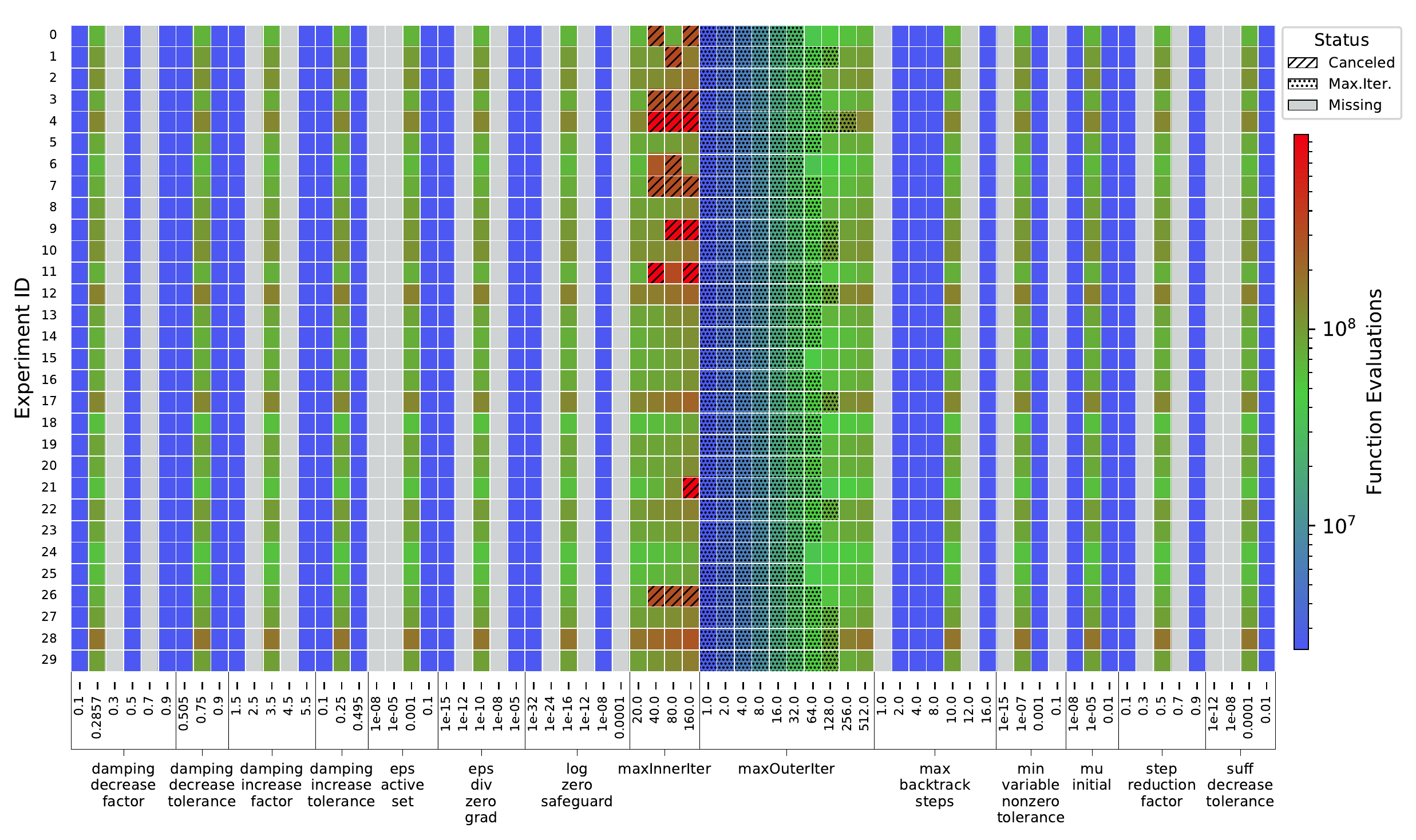}
		\subcaption{{\it PDNR}, {\it nell-2}}
		\label{fig:results:heatmap:arm:pdnr:nell}
	\end{subfigure}
	~
	\begin{subfigure}{\textwidth}
		\includegraphics[width=\textwidth]{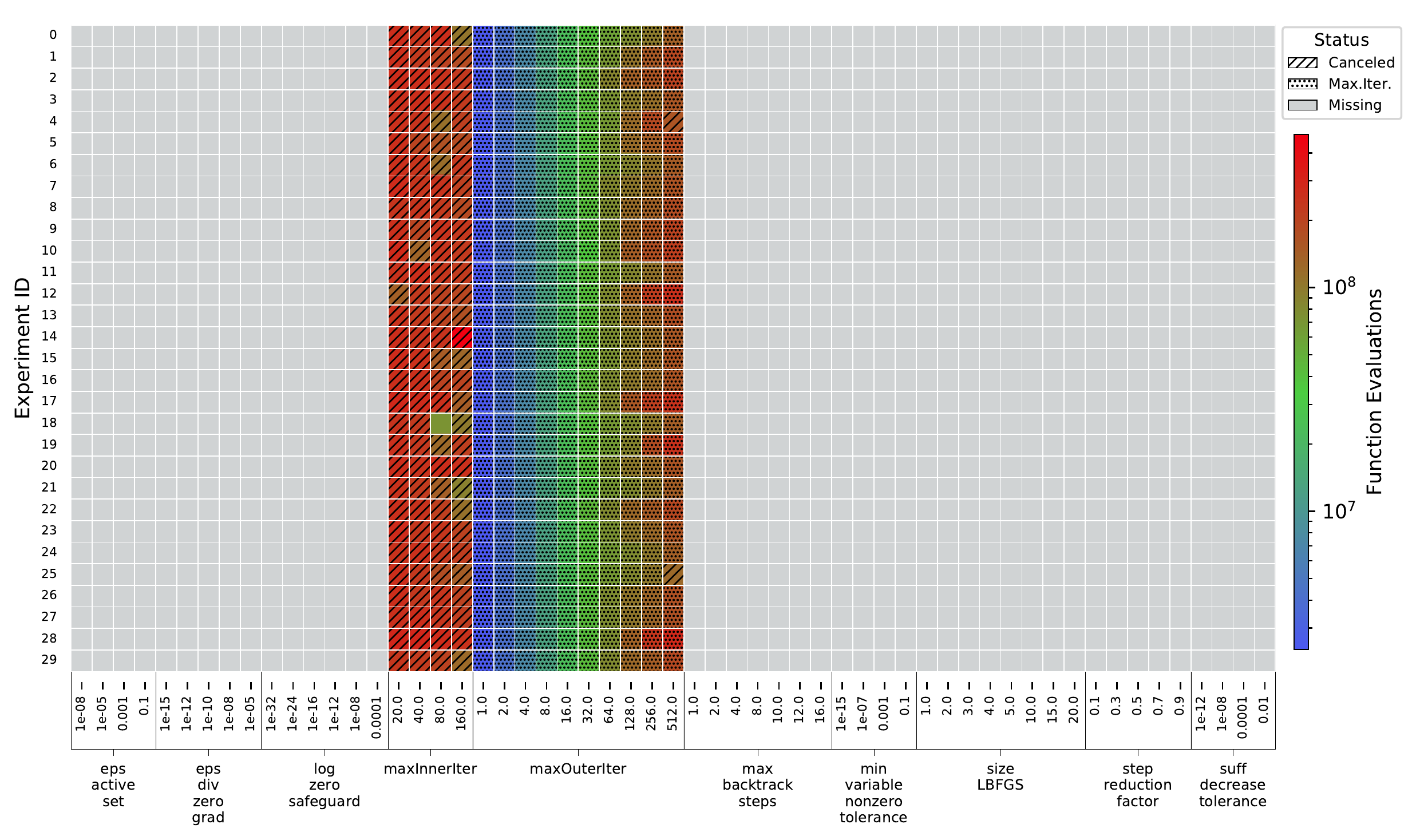}
		\subcaption{{\it PQNR}, {\it nell-2}}
		\label{fig:results:heatmap:arm:pqnr:nell}
	\end{subfigure}
\end{figure}
\clearpage

\begin{figure}[!htp]
	\centering
	\caption{Experiment outcomes for {\it chicago-crime-comm} data on IBM platform.}
	\label{fig:results:heatmap:ibm:chicago}
	\begin{subfigure}{\textwidth}
		\includegraphics[width=\textwidth]{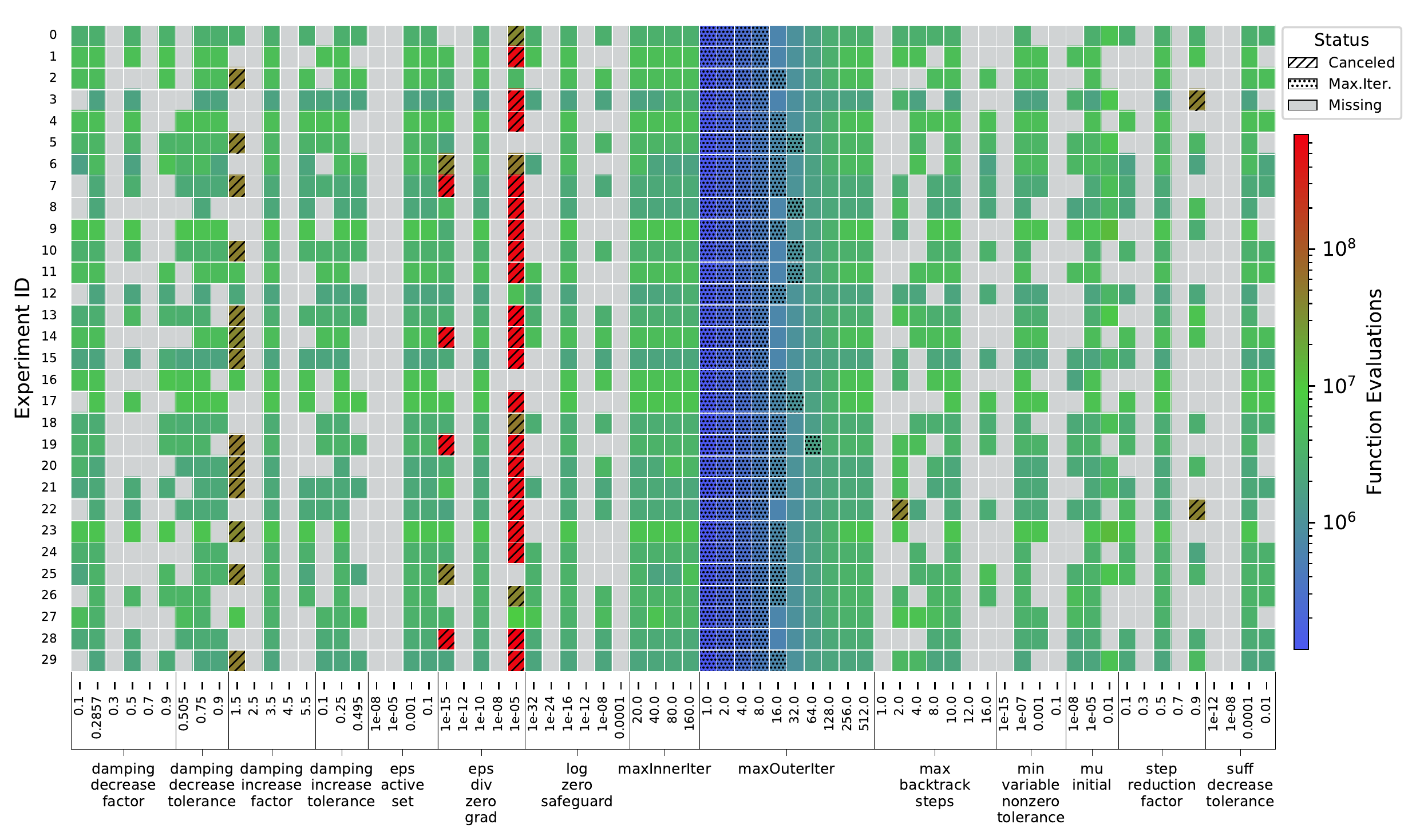}
		\subcaption{{\it PDNR}, {\it chicago-crime-comm}}
		\label{fig:results:heatmap:ibm:pdnr:chicago}
	\end{subfigure}
	~
	\begin{subfigure}{\textwidth}
		\includegraphics[width=\textwidth]{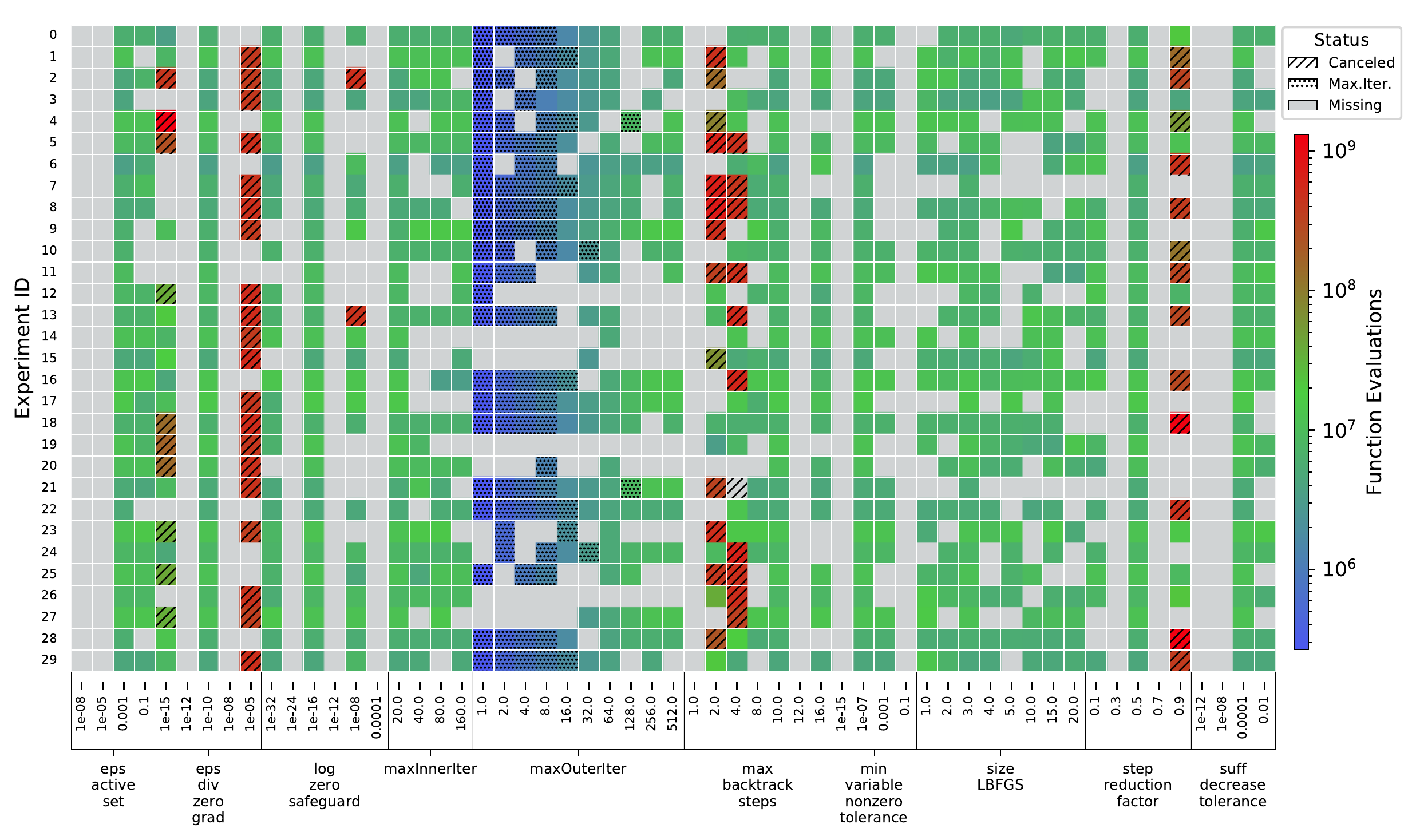}
		\subcaption{{\it PQNR}, {\it chicago-crime-comm}}
		\label{fig:results:heatmap:ibm:pqnr:chicago}
	\end{subfigure}
\end{figure}
\clearpage
\begin{figure}[!htp]
	\centering
	\caption{Experiment outcomes for {\it lbnl-network} data on IBM platform.}
	\label{fig:results:heatmap:ibm:lbnl}
	\begin{subfigure}{\textwidth}
		\includegraphics[width=\textwidth]{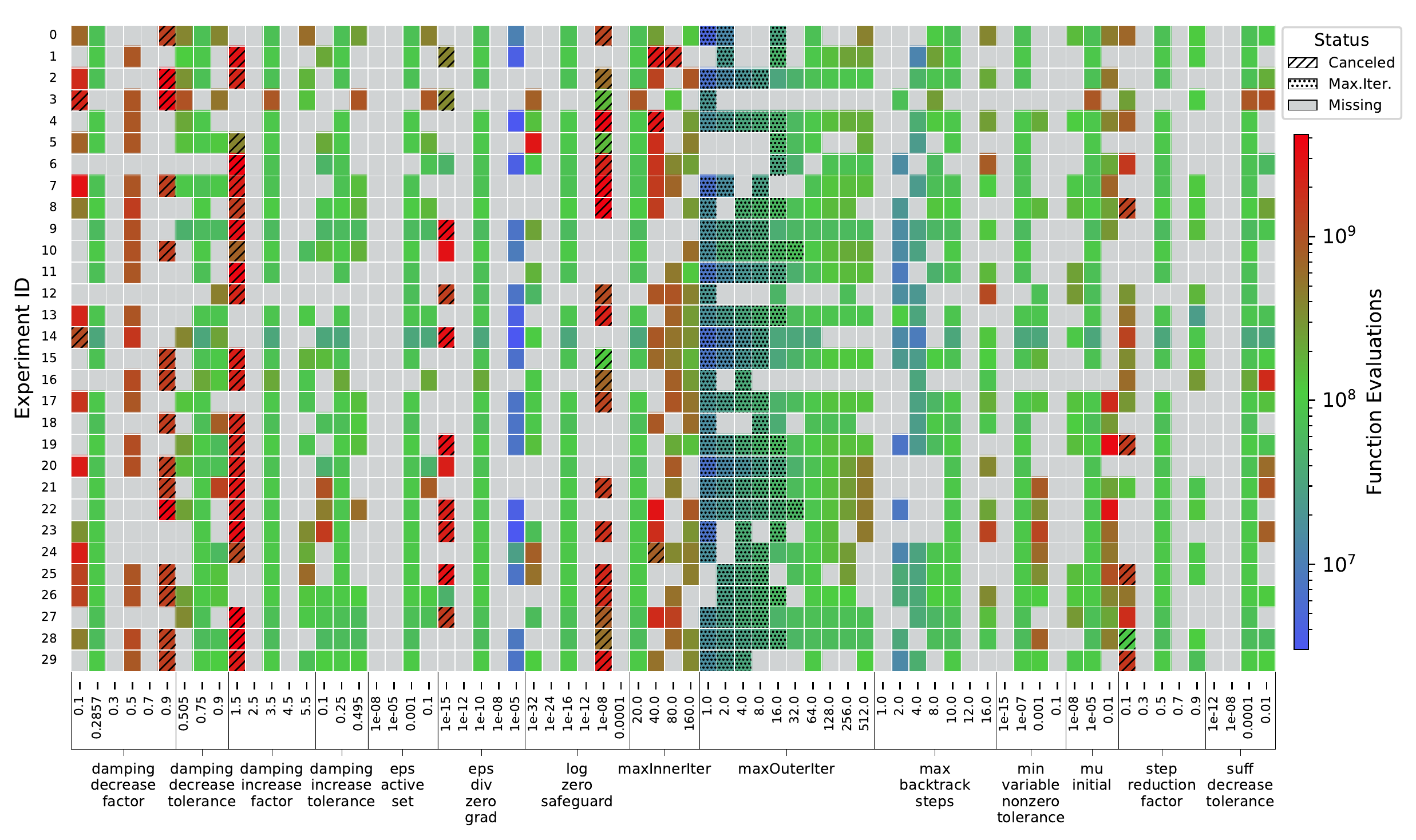}
		\subcaption{{\it PDNR}, {\it lbnl-network}}
		\label{fig:results:heatmap:ibm:pdnr:lbnl}
	\end{subfigure}
	~
	\begin{subfigure}{\textwidth}
		\includegraphics[width=\textwidth]{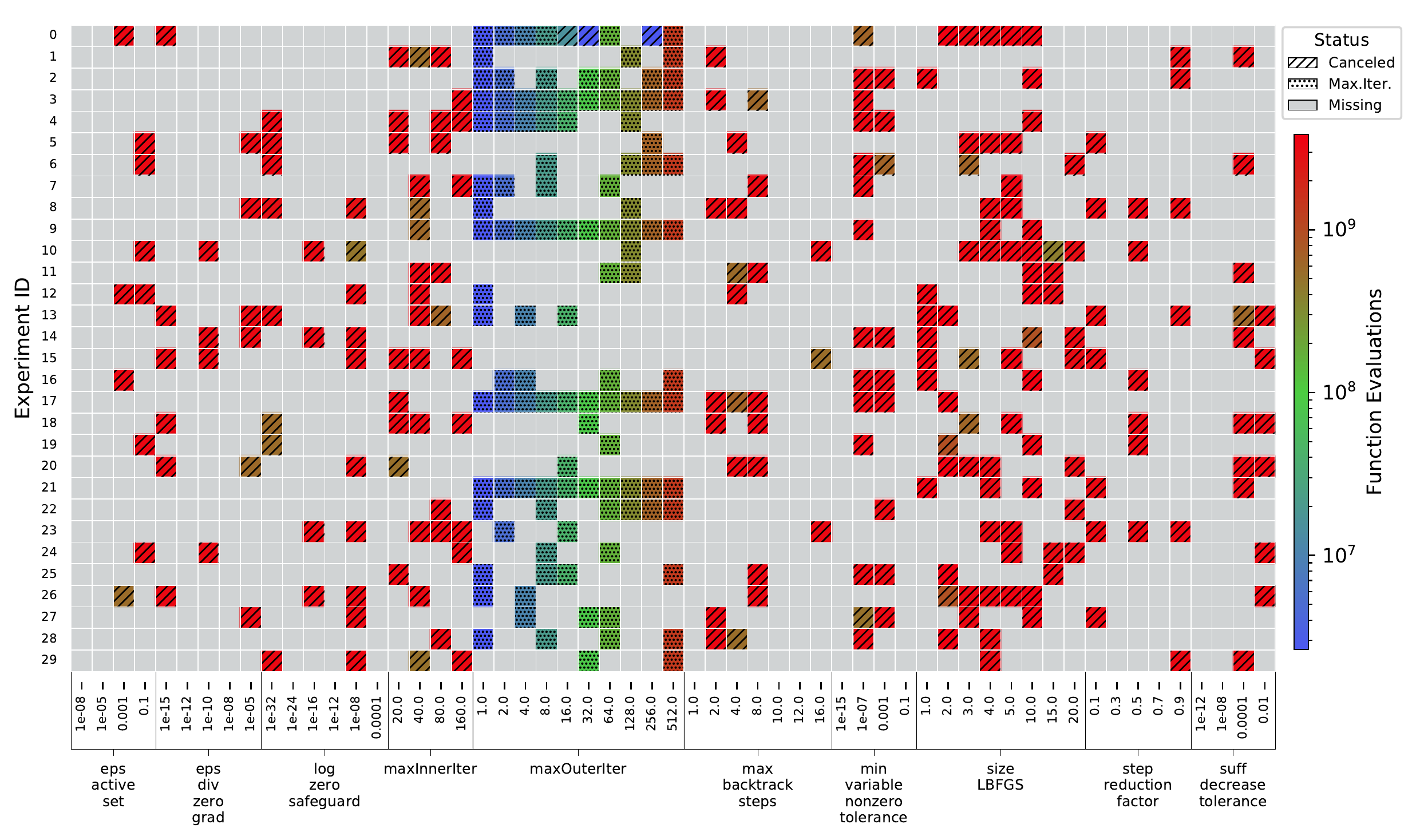}
		\subcaption{{\it PQNR}, {\it lbnl-network}}
		\label{fig:results:heatmap:ibm:pqnr:lbnl}
	\end{subfigure}
\end{figure}
\clearpage
\begin{figure}[!htp]
	\centering
	\caption{Experiment outcomes for {\it nell-2} data on IBM platform.}
	\label{fig:results:heatmap:ibm:nell}
	\begin{subfigure}{\textwidth}
		\includegraphics[width=\textwidth]{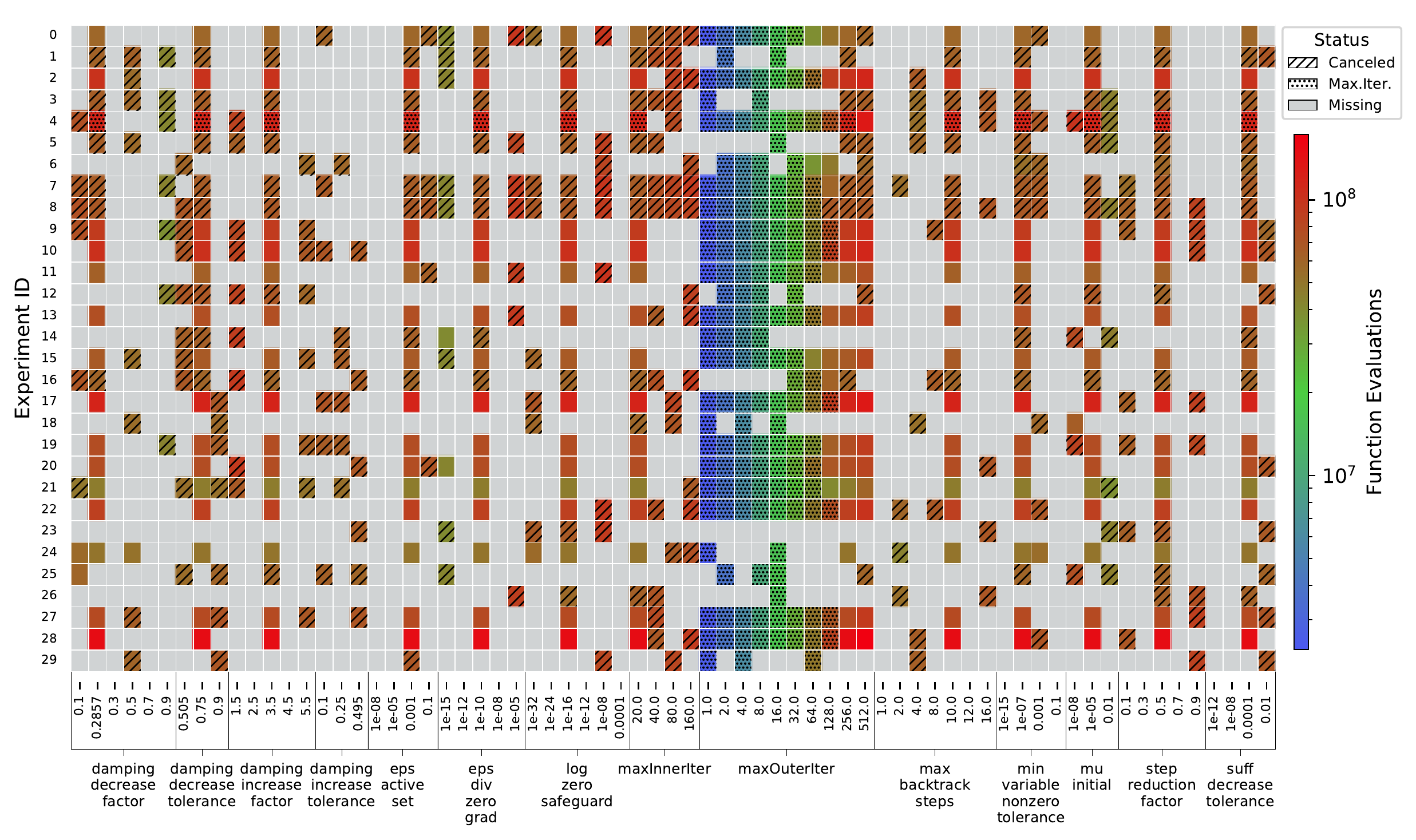}
		\subcaption{{\it PDNR}, {\it nell-2}}
		\label{fig:results:heatmap:ibm:pdnr:nell}
	\end{subfigure}
	~
	\begin{subfigure}{\textwidth}
		\includegraphics[width=\textwidth]{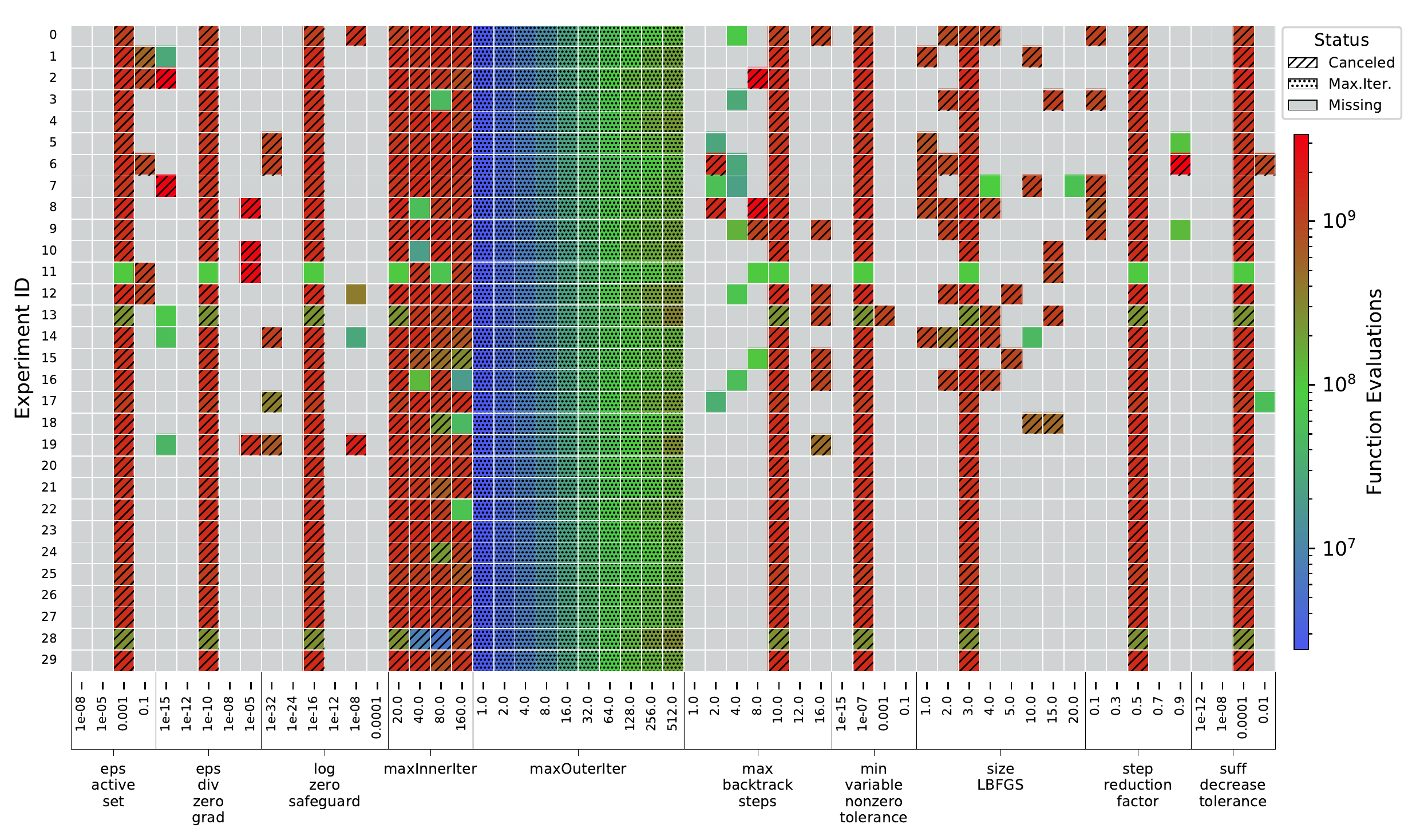}
		\subcaption{{\it PQNR}, {\it nell-2}}
		\label{fig:results:heatmap:ibm:pqnr:nell}
	\end{subfigure}
\end{figure}
\clearpage

\subsection{Average Convergence Results}\label{app:supplementary:convergence}
Below are tables presenting the number of of objective function evaluations for 
PDNR and PQNR to reach convergence, averaged across hardware platforms and 
reported per data set. $N$ denotes number of {\it converged} experiments per 
set of experiments. Mean and 95\% confidence intervals (CI) are given for 
objection function evaluation results.

\begin{table}[!ht]
	\centering
	\caption{Average convergence behavior for the {\it chicago-crime-comm} data 
		set. Confidence intervals are percentages of the mean.}
	\label{tab:results:chicago-crime-comm:fevals}
	{\resizebox{0.47\textwidth}{!}{%
			\begin{tabular}[t]{lrrrrr}
				\toprule
				Parameter&Value&Solver&$N$&\multicolumn{2}{c}{Function 
				Evaluations} \\
				&&&&Mean & 95\% CI \\\midrule
				\multirow{3}*{\texttt{mu\_initial}}
				&\multirow{1}*{1e-08}&PDNR&78&3.30e+06&$\pm 8.2$\\
				&\multirow{1}*{1e-05}&PDNR&90&3.54e+06&$\pm 7.9$\\
				&\multirow{1}*{0.01}&PDNR&76&6.80e+06&$\pm 10.2$\\
				\midrule
				\multirow{6}*{\texttt{damping\_decrease\_factor}}
				&\multirow{1}*{0.1}&PDNR&82&3.49e+06&$\pm 8.3$\\
				&\multirow{1}*{0.2857}&PDNR&90&3.54e+06&$\pm 7.9$\\
				&\multirow{1}*{0.3}&PDNR&27&3.49e+06&$\pm 16.9$\\
				&\multirow{1}*{0.5}&PDNR&79&3.55e+06&$\pm 8.8$\\
				&\multirow{1}*{0.7}&PDNR&26&3.80e+06&$\pm 15.4$\\
				&\multirow{1}*{0.9}&PDNR&77&4.02e+06&$\pm 8.1$\\
				\midrule
				\multirow{3}*{\texttt{damping\_decrease\_tolerance}}
				&\multirow{1}*{0.505}&PDNR&77&3.56e+06&$\pm 8.6$\\
				&\multirow{1}*{0.75}&PDNR&90&3.54e+06&$\pm 7.9$\\
				&\multirow{1}*{0.9}&PDNR&82&3.45e+06&$\pm 8.6$\\
				\midrule
				\multirow{5}*{\texttt{damping\_increase\_factor}}
				&\multirow{1}*{1.5}&PDNR&31&3.97e+06&$\pm 13.8$\\
				&\multirow{1}*{2.5}&PDNR&30&3.53e+06&$\pm 14.8$\\
				&\multirow{1}*{3.5}&PDNR&90&3.54e+06&$\pm 7.9$\\
				&\multirow{1}*{4.5}&PDNR&30&3.45e+06&$\pm 14.7$\\
				&\multirow{1}*{5.5}&PDNR&80&3.42e+06&$\pm 9.1$\\
				\midrule
				\multirow{3}*{\texttt{damping\_increase\_tolerance}}
				&\multirow{1}*{0.1}&PDNR&78&3.46e+06&$\pm 8.5$\\
				&\multirow{1}*{0.25}&PDNR&90&3.54e+06&$\pm 7.9$\\
				&\multirow{1}*{0.495}&PDNR&76&3.47e+06&$\pm 9.0$\\
				\midrule
				\multirow{9}*{\texttt{eps\_active\_set}}
				&\multirow{2}*{1e-08}&PDNR&30&3.53e+06&$\pm 14.3$\\
				&&PQNR&30&8.73e+06&$\pm 14.5$\\
				\cmidrule{2-6}
				&\multirow{2}*{1e-05}&PDNR&30&3.53e+06&$\pm 14.3$\\
				&&PQNR&30&8.28e+06&$\pm 14.5$\\
				\cmidrule{2-6}
				&\multirow{2}*{0.001}&PDNR&90&3.54e+06&$\pm 7.9$\\
				&&PQNR&88&8.37e+06&$\pm 8.0$\\
				\cmidrule{2-6}
				&\multirow{2}*{0.1}&PDNR&90&3.54e+06&$\pm 7.9$\\
				&&PQNR&53&8.16e+06&$\pm 10.6$\\
				\midrule
				\multirow{10.5}*{\texttt{eps\_div\_zero\_grad}}
				&\multirow{2}*{1e-15}
				&PDNR&67&3.20e+06&$\pm 9.3$\\
				&&PQNR&20&1.14e+07&$\pm 24.0$\\
				\cmidrule{2-6}
				&\multirow{2}*{1e-12}&PDNR&30&3.59e+06&$\pm 13.5$\\
				&&PQNR&25&8.28e+06&$\pm 15.3$\\
				\cmidrule{2-6}
				&\multirow{2}*{1e-10}&PDNR&90&3.54e+06&$\pm 7.9$\\
				&&PQNR&88&8.37e+06&$\pm 8.0$\\
				\cmidrule{2-6}
				&\multirow{2}*{1e-08}&PDNR&30&4.57e+06&$\pm 14.9$\\
				&&PQNR&30&8.45e+06&$\pm 14.0$\\
				\cmidrule{2-6}
				&\multirow{1}*{1e-05}&PDNR&11&5.37e+06&$\pm 26.8$\\
				\midrule
				\multirow{12.5}*{\texttt{log\_zero\_safeguard}}
				&\multirow{2}*{1e-32}&PDNR&76&3.52e+06&$\pm 8.9$\\
				&&PQNR&56&7.90e+06&$\pm 10.6$\\
				\cmidrule{2-6}
				&\multirow{2}*{1e-24}&PDNR&30&3.55e+06&$\pm 15.1$\\
				&&PQNR&30&8.49e+06&$\pm 14.5$\\
				\cmidrule{2-6}
				&\multirow{2}*{1e-16}&PDNR&90&3.54e+06&$\pm 7.9$\\
				&&PQNR&88&8.37e+06&$\pm 8.0$\\
				\cmidrule{2-6}
				&\multirow{2}*{1e-12}&PDNR&30&3.45e+06&$\pm 15.3$\\
				&&PQNR&30&8.01e+06&$\pm 15.0$\\
				\cmidrule{2-6}
				&\multirow{2}*{1e-08}&PDNR&78&3.46e+06&$\pm 8.5$\\
				&&PQNR&47&8.81e+06&$\pm 11.2$\\
				\cmidrule{2-6}
				&\multirow{1}*{0.0001}&PQNR&8&6.48e+07&$\pm 73.8$\\
				\midrule
				\multirow{9}*{\texttt{max\_inner\_iterations}}
				&\multirow{2}*{20}&PDNR&90&3.54e+06&$\pm 7.9$\\
				&&PQNR&88&8.37e+06&$\pm 8.0$\\
				\cmidrule{2-6}
				&\multirow{2}*{40}&PDNR&90&3.51e+06&$\pm 8.7$\\
				&&PQNR&49&8.42e+06&$\pm 11.7$\\
				\cmidrule{2-6}
				&\multirow{2}*{80}&PDNR&90&3.59e+06&$\pm 8.1$\\
				&&PQNR&52&8.68e+06&$\pm 11.3$\\
				\cmidrule{2-6}
				&\multirow{2}*{160}&PDNR&90&3.68e+06&$\pm 7.6$\\
				&&PQNR&52&8.26e+06&$\pm 10.8$\\
				\bottomrule
	\end{tabular}}}%
	{\resizebox{0.47\textwidth}{!}{%
			\begin{tabular}[t]{lrrrrr}
				\toprule
				Parameter&Value&Solver&$N$&\multicolumn{2}{c}{Function 
				Evaluations} \\
				&&&&Mean & 95\% CI \\[0.15em]\midrule
				\multirow{14.5}*{\texttt{max\_outer\_iterations}}
				&\multirow{1}*{8}&PQNR&3&1.03e+06&$\pm 26.9$\\
				\cmidrule{2-6}
				&\multirow{2}*{16}&PDNR&36&6.43e+05&$\pm 1.4$\\
				&&PQNR&48&1.94e+06&$\pm 3.9$\\
				\cmidrule{2-6}
				&\multirow{2}*{32}&PDNR&75&1.03e+06&$\pm 1.2$\\
				&&PQNR&74&3.06e+06&$\pm 2.5$\\
				\cmidrule{2-6}
				&\multirow{2}*{64}&PDNR&87&1.69e+06&$\pm 1.0$\\
				&&PQNR&83&4.66e+06&$\pm 2.7$\\
				\cmidrule{2-6}
				&\multirow{2}*{128}&PDNR&90&2.69e+06&$\pm 3.3$\\
				&&PQNR&66&6.93e+06&$\pm 5.8$\\
				\cmidrule{2-6}
				&\multirow{2}*{256}&PDNR&90&3.48e+06&$\pm 7.5$\\
				&&PQNR&73&8.50e+06&$\pm 9.0$\\
				\cmidrule{2-6}
				&\multirow{2}*{512}&PDNR&90&3.54e+06&$\pm 7.9$\\
				&&PQNR&77&8.51e+06&$\pm 8.6$\\
				\midrule
				\multirow{16}*{\texttt{max\_backtrack\_steps}}
				&\multirow{2}*{1}&PDNR&29&3.64e+06&$\pm 12.6$\\
				&&PQNR&7&1.52e+07&$\pm 46.4$\\
				\cmidrule{2-6}
				&\multirow{2}*{2}&PDNR&75&3.87e+06&$\pm 7.5$\\
				&&PQNR&14&1.69e+07&$\pm 40.9$\\
				\cmidrule{2-6}
				&\multirow{2}*{4}&PDNR&77&3.85e+06&$\pm 7.6$\\
				&&PQNR&29&1.04e+07&$\pm 23.0$\\
				\cmidrule{2-6}
				&\multirow{2}*{8}&PDNR&80&3.51e+06&$\pm 8.4$\\
				&&PQNR&48&8.06e+06&$\pm 12.1$\\
				\cmidrule{2-6}
				&\multirow{2}*{10}&PDNR&90&3.54e+06&$\pm 7.9$\\
				&&PQNR&88&8.37e+06&$\pm 8.0$\\
				\cmidrule{2-6}
				&\multirow{2}*{12}&PDNR&30&3.51e+06&$\pm 15.8$\\
				&&PQNR&29&8.08e+06&$\pm 12.0$\\
				\cmidrule{2-6}
				&\multirow{2}*{16}&PDNR&75&3.57e+06&$\pm 9.1$\\
				&&PQNR&54&8.44e+06&$\pm 10.5$\\
				\midrule
				\multirow{9}*{\texttt{min\_variable\_nonzero\_tolerance}}
				&\multirow{2}*{1e-15}&PDNR&30&3.54e+06&$\pm 14.2$\\
				&&PQNR&30&8.73e+06&$\pm 14.5$\\
				\cmidrule{2-6}
				&\multirow{2}*{1e-07}&PDNR&90&3.54e+06&$\pm 7.9$\\
				&&PQNR&88&8.37e+06&$\pm 8.0$\\
				\cmidrule{2-6}
				&\multirow{2}*{0.001}&PDNR&80&3.56e+06&$\pm 8.5$\\
				&&PQNR&47&8.26e+06&$\pm 11.1$\\
				\cmidrule{2-6}
				&\multirow{2}*{0.1}&PDNR&30&3.54e+06&$\pm 14.2$\\
				&&PQNR&30&8.46e+06&$\pm 15.3$\\
				\midrule
				\multirow{8}*{\texttt{size\_LBFGS}}
				&\multirow{1}*{1}&PQNR&51&9.30e+06&$\pm 11.4$\\
				&\multirow{1}*{2}&PQNR&51&7.89e+06&$\pm 10.9$\\
				&\multirow{1}*{3}&PQNR&88&8.37e+06&$\pm 8.0$\\
				&\multirow{1}*{4}&PQNR&50&7.44e+06&$\pm 10.0$\\
				&\multirow{1}*{5}&PQNR&49&8.14e+06&$\pm 11.8$\\
				&\multirow{1}*{10}&PQNR&47&7.26e+06&$\pm 11.2$\\
				&\multirow{1}*{10}&PQNR&55&8.12e+06&$\pm 10.2$\\
				&\multirow{1}*{20}&PQNR&49&7.23e+06&$\pm 11.8$\\
				\midrule
				\multirow{11}*{\texttt{step\_reduction\_factor}}
				&\multirow{2}*{0.1}&PDNR&79&3.61e+06&$\pm 8.2$\\
				&&PQNR&51&8.42e+06&$\pm 10.6$\\
				\cmidrule{2-6}
				&\multirow{2}*{0.3}&PDNR&30&3.59e+06&$\pm 13.4$\\
				&&PQNR&30&8.11e+06&$\pm 15.4$\\
				\cmidrule{2-6}
				&\multirow{2}*{0.5}&PDNR&90&3.54e+06&$\pm 7.9$\\
				&&PQNR&88&8.37e+06&$\pm 8.0$\\
				\cmidrule{2-6}
				&\multirow{2}*{0.7}&PDNR&30&3.88e+06&$\pm 14.6$\\
				&&PQNR&26&7.97e+06&$\pm 18.3$\\
				\cmidrule{2-6}
				&\multirow{2}*{0.9}&PDNR&71&3.96e+06&$\pm 8.4$\\
				&&PQNR&12&1.20e+07&$\pm 29.8$\\
				\midrule
				\multirow{9}*{\texttt{suff\_decrease\_tolerance}}
				&\multirow{2}*{1e-12}&PDNR&30&3.52e+06&$\pm 14.6$\\
				&&PQNR&30&8.19e+06&$\pm 14.5$\\
				\cmidrule{2-6}
				&\multirow{2}*{1e-08}&PDNR&30&3.52e+06&$\pm 14.6$\\
				&&PQNR&30&8.13e+06&$\pm 14.7$\\
				\cmidrule{2-6}
				&\multirow{2}*{0.0001}&PDNR&90&3.54e+06&$\pm 7.9$\\
				&&PQNR&88&8.37e+06&$\pm 8.0$\\
				\cmidrule{2-6}
				&\multirow{2}*{0.01}&PDNR&78&3.50e+06&$\pm 8.7$\\
				&&PQNR&51&8.38e+06&$\pm 10.9$\\
				\bottomrule
	\end{tabular}}}%
\end{table}
\begin{table}[!ht]
	\centering
	\caption{Average convergence behavior for the \textit{lbnl-network} data 
	set. 
		Confidence intervals are percentages of the mean.}
	\label{tab:results:lbnl-network:fevals}
	{\resizebox{0.48\textwidth}{!}{%
			\begin{tabular}[t]{lrrrrr}
				\toprule
				Parameter&Value&Solver&$N$&\multicolumn{2}{c}{Function 
				Evaluations} \\
				&&&&Mean & 95\% CI \\\midrule
				\multirow{3}*{\texttt{mu\_initial}}
				&\multirow{1}*{1e-08}&PDNR&74&1.84e+08&$\pm 23.0$\\
				&\multirow{1}*{1e-05}&PDNR&89&2.34e+08&$\pm 25.1$\\
				&\multirow{1}*{0.01}&PDNR&78&7.83e+08&$\pm 23.9$\\
				\midrule
				\multirow{4}*{\texttt{damping\_decrease\_factor}}
				&\multirow{1}*{0.1}&PDNR&71&1.56e+09&$\pm 14.4$\\
				&\multirow{1}*{0.2857}&PDNR&86&2.30e+08&$\pm 25.4$\\
				&\multirow{1}*{0.3}&PDNR&20&3.29e+08&$\pm 28.6$\\
				&\multirow{1}*{0.5}&PDNR&79&1.02e+09&$\pm 4.2$\\
				\midrule
				\multirow{3}*{\texttt{damping\_decrease\_tolerance}}
				&\multirow{1}*{0.505}&PDNR&75&2.70e+08&$\pm 21.3$\\
				&\multirow{1}*{0.75}&PDNR&87&2.30e+08&$\pm 25.2$\\
				&\multirow{1}*{0.9}&PDNR&79&2.63e+08&$\pm 25.4$\\
				\midrule
				\multirow{4}*{\texttt{damping\_increase\_factor}}
				&\multirow{1}*{2.5}&PDNR&30&1.58e+08&$\pm 32.1$\\
				&\multirow{1}*{3.5}&PDNR&89&2.36e+08&$\pm 24.9$\\
				&\multirow{1}*{4.5}&PDNR&30&1.85e+08&$\pm 41.5$\\
				&\multirow{1}*{5.5}&PDNR&75&1.73e+08&$\pm 22.1$\\
				\midrule
				\multirow{3}*{\texttt{damping\_increase\_tolerance}}
				&\multirow{1}*{0.1}&PDNR&77&2.60e+08&$\pm 30.5$\\
				&\multirow{1}*{0.25}&PDNR&88&2.28e+08&$\pm 25.1$\\
				&\multirow{1}*{0.495}&PDNR&74&2.64e+08&$\pm 27.9$\\
				\midrule
				\multirow{4}*{\texttt{eps\_active\_set}}
				&\multirow{1}*{1e-08}&PDNR&30&2.80e+08&$\pm 39.8$\\
				&\multirow{1}*{1e-05}&PDNR&30&2.78e+08&$\pm 39.7$\\
				&\multirow{1}*{0.001}&PDNR&87&2.28e+08&$\pm 25.4$\\
				&\multirow{1}*{0.1}&PDNR&73&2.49e+08&$\pm 23.9$\\
				\midrule
				\multirow{5}*{\texttt{eps\_div\_zero\_grad}}
				&\multirow{1}*{1e-15}&PDNR&33&1.85e+09&$\pm 29.8$\\
				&\multirow{1}*{1e-12}&PDNR&29&1.38e+09&$\pm 21.2$\\
				&\multirow{1}*{1e-10}&PDNR&89&2.27e+08&$\pm 25.1$\\
				&\multirow{1}*{1e-08}&PDNR&30&2.45e+07&$\pm 22.7$\\
				&\multirow{1}*{1e-05}&PDNR&79&7.28e+06&$\pm 14.3$\\
				\midrule
				\multirow{4}*{\texttt{log\_zero\_safeguard}}
				&\multirow{1}*{1e-32}&PDNR&77&2.95e+08&$\pm 43.7$\\
				&\multirow{1}*{1e-24}&PDNR&30&2.82e+08&$\pm 40.2$\\
				&\multirow{1}*{1e-16}&PDNR&87&2.29e+08&$\pm 25.4$\\
				&\multirow{1}*{1e-12}&PDNR&30&4.20e+08&$\pm 66.9$\\
				\bottomrule
	\end{tabular}}}%
	{\resizebox{0.48\textwidth}{!}{%
			\begin{tabular}[t]{lrrrrr}
				\toprule
				Parameter&Value&Solver&$N$&\multicolumn{2}{c}{Function 
				Evaluations} \\
				&&&&Mean & 95\% CI \\[0.15em]\midrule
				\multirow{4}*{\texttt{max\_inner\_iterations}}
				&\multirow{1}*{20}&PDNR&88&2.36e+08&$\pm 25.1$\\
				&\multirow{1}*{40}&PDNR&69&1.13e+09&$\pm 13.6$\\
				&\multirow{1}*{80}&PDNR&77&7.31e+08&$\pm 13.9$\\
				&\multirow{1}*{160}&PDNR&82&4.00e+08&$\pm 13.1$\\
				\midrule
				\multirow{6}*{\texttt{max\_outer\_iterations}}
				&\multirow{1}*{16}&PDNR&35&4.32e+07&$\pm 6.4$\\
				&\multirow{1}*{32}&PDNR&69&5.85e+07&$\pm 4.2$\\
				&\multirow{1}*{64}&PDNR&86&8.33e+07&$\pm 3.9$\\
				&\multirow{1}*{128}&PDNR&82&1.21e+08&$\pm 6.5$\\
				&\multirow{1}*{256}&PDNR&83&1.62e+08&$\pm 10.4$\\
				&\multirow{1}*{512}&PDNR&81&2.20e+08&$\pm 16.1$\\
				\midrule
				\multirow{7}*{\texttt{max\_backtrack\_steps}}
				&\multirow{1}*{1}&PDNR&25&1.05e+08&$\pm 15.7$\\
				&\multirow{1}*{2}&PDNR&77&2.85e+07&$\pm 19.8$\\
				&\multirow{1}*{4}&PDNR&80&3.22e+07&$\pm 10.6$\\
				&\multirow{1}*{8}&PDNR&79&1.39e+08&$\pm 21.3$\\
				&\multirow{1}*{10}&PDNR&86&2.30e+08&$\pm 25.4$\\
				&\multirow{1}*{12}&PDNR&28&3.76e+08&$\pm 56.5$\\
				&\multirow{1}*{16}&PDNR&79&3.13e+08&$\pm 23.8$\\
				\midrule
				\multirow{4}*{\texttt{min\_variable\_nonzero\_tolerance}}
				&\multirow{1}*{1e-15}&PDNR&30&2.88e+08&$\pm 39.9$\\
				&\multirow{1}*{1e-07}&PDNR&87&2.29e+08&$\pm 25.4$\\
				&\multirow{1}*{0.001}&PDNR&75&3.12e+08&$\pm 23.6$\\
				&\multirow{1}*{0.1}&PDNR&25&2.66e+08&$\pm 43.3$\\
				\midrule
				\multirow{5}*{\texttt{step\_reduction\_factor}}
				&\multirow{1}*{0.1}&PDNR&59&1.01e+09&$\pm 23.4$\\
				&\multirow{1}*{0.3}&PDNR&29&3.61e+08&$\pm 39.5$\\
				&\multirow{1}*{0.5}&PDNR&87&2.29e+08&$\pm 25.4$\\
				&\multirow{1}*{0.7}&PDNR&29&7.52e+07&$\pm 30.1$\\
				&\multirow{1}*{0.9}&PDNR&77&1.65e+08&$\pm 15.2$\\
				\midrule
				\multirow{4}*{\texttt{suff\_decrease\_tolerance}}
				&\multirow{1}*{1e-12}&PDNR&30&2.98e+08&$\pm 40.7$\\
				&\multirow{1}*{1e-08}&PDNR&30&2.91e+08&$\pm 44.0$\\
				&\multirow{1}*{0.0001}&PDNR&90&2.34e+08&$\pm 24.8$\\
				&\multirow{1}*{0.01}&PDNR&78&3.78e+08&$\pm 29.8$\\
				\bottomrule
	\end{tabular}}}%
\end{table}
\begin{table}[!ht]
	\centering
	\caption{Average convergence behavior for the \textit{nell-2} data set. 
		Confidence intervals are percentages of the mean.}
	\label{tab:results:nell-2:fevals}
	{\resizebox{0.48\textwidth}{!}{%
			\begin{tabular}[t]{lrrrrr}
				\toprule
				Parameter&Value&Solver&$N$&\multicolumn{2}{c}{Function 
				Evaluations} \\
				&&&&Mean & 95\% CI \\\midrule
				\multirow{2}*{\texttt{mu\_initial}}
				&\multirow{1}*{1e-08}&PDNR&14&1.01e+08&$\pm 19.9$\\
				&\multirow{1}*{1e-05}&PDNR&65&9.08e+07&$\pm 7.1$\\
				\midrule
				\multirow{6}*{\texttt{damping\_decrease\_factor}}
				&\multirow{1}*{0.1}&PDNR&20&9.74e+07&$\pm 15.4$\\
				&\multirow{1}*{0.2857}&PDNR&74&8.94e+07&$\pm 6.8$\\
				&\multirow{1}*{0.3}&PDNR&11&1.00e+08&$\pm 20.3$\\
				&\multirow{1}*{0.5}&PDNR&24&9.34e+07&$\pm 13.7$\\
				&\multirow{1}*{0.7}&PDNR&16&8.62e+07&$\pm 17.6$\\
				&\multirow{1}*{0.9}&PDNR&19&8.68e+07&$\pm 15.9$\\
				\midrule
				\multirow{3}*{\texttt{damping\_decrease\_tolerance}}
				&\multirow{1}*{0.505}&PDNR&18&9.41e+07&$\pm 13.6$\\
				&\multirow{1}*{0.75}&PDNR&64&9.06e+07&$\pm 7.2$\\
				&\multirow{1}*{0.9}&PDNR&21&9.76e+07&$\pm 14.5$\\
				\midrule
				\multirow{4}*{\texttt{damping\_increase\_factor}}
				&\multirow{1}*{2.5}&PDNR&15&9.31e+07&$\pm 10.8$\\
				&\multirow{1}*{3.5}&PDNR&64&9.06e+07&$\pm 7.2$\\
				&\multirow{1}*{4.5}&PDNR&17&9.49e+07&$\pm 12.8$\\
				&\multirow{1}*{5.5}&PDNR&21&9.06e+07&$\pm 15.0$\\
				\midrule
				\multirow{3}*{\texttt{damping\_increase\_tolerance}}
				&\multirow{1}*{0.1}&PDNR&19&1.04e+08&$\pm 14.1$\\
				&\multirow{1}*{0.25}&PDNR&49&9.29e+07&$\pm 7.9$\\
				&\multirow{1}*{0.49}&PDNR&16&9.03e+07&$\pm 13.2$\\
				\midrule
				\multirow{5}*{\texttt{eps\_active\_set}}
				&\multirow{1}*{1e-08}&PDNR&20&9.18e+07&$\pm 11.5$\\
				&\multirow{1}*{1e-05}&PDNR&20&9.55e+07&$\pm 12.3$\\
				\cmidrule{2-6}
				&\multirow{2}*{0.001}&PDNR&74&8.94e+07&$\pm 6.8$\\
				&&PQNR&1&9.69e+07&$\pm 0.0$\\
				\cmidrule{2-6}
				&\multirow{1}*{0.1}&PDNR&20&9.16e+07&$\pm 11.1$\\
				\midrule
				\multirow{7}*{\texttt{eps\_div\_zero\_grad}}
				&\multirow{2}*{1e-15}&PDNR&17&5.37e+07&$\pm 13.9$\\
				&&PQNR&21&1.86e+08&$\pm 39.8$\\
				\cmidrule{2-6}
				&\multirow{2}*{1e-12}&PDNR&30&6.96e+07&$\pm 9.6$\\
				&&PQNR&9&1.24e+08&$\pm 40.1$\\
				\cmidrule{2-6}
				&\multirow{2}*{1e-10}&PDNR&74&8.94e+07&$\pm 6.8$\\
				&&PQNR&1&9.69e+07&$\pm 0.0$\\
				\midrule
				\multirow{9}*{\texttt{log\_zero\_safeguard}}
				&\multirow{1}*{1e-32}&PDNR&20&9.94e+07&$\pm 14.8$\\
				&\multirow{1}*{1e-24}&PDNR&16&9.55e+07&$\pm 12.5$\\
				\cmidrule{2-6}
				&\multirow{2}*{1e-16}&PDNR&74&8.94e+07&$\pm 6.8$\\
				&&PQNR&1&9.69e+07&$\pm 0.0$\\
				\cmidrule{2-6}
				&\multirow{2}*{1e-12}&PDNR&19&1.07e+08&$\pm 16.1$\\
				&&PQNR&1&9.56e+07&$\pm 0.0$\\
				\cmidrule{2-6}
				&\multirow{1}*{1e-08}&PQNR&14&2.62e+08&$\pm 30.4$\\
				&\multirow{1}*{0.0001}&PQNR&21&3.03e+08&$\pm 37.1$\\
				\bottomrule
	\end{tabular}}}%
	{\resizebox{0.48\textwidth}{!}{%
			\begin{tabular}[t]{lrrrrr}
				\toprule
				Parameter&Value&Solver&$N$&\multicolumn{2}{c}{Function 
				Evaluations} \\
				&&&&Mean & 95\% CI \\[0.15em]\midrule
				\multirow{9}*{\texttt{max\_inner\_iterations}}
				&\multirow{2}*{20}&PDNR&74&8.94e+07&$\pm 6.8$\\
				&&PQNR&1&9.69e+07&$\pm 0.0$\\
				\cmidrule{2-6}
				&\multirow{2}*{40}&PDNR&48&1.19e+08&$\pm 16.5$\\
				&&PQNR&3&6.77e+07&$\pm 199.0$\\
				\cmidrule{2-6}
				&\multirow{2}*{80}&PDNR&46&1.28e+08&$\pm 13.8$\\
				&&PQNR&4&6.54e+07&$\pm 28.7$\\
				\cmidrule{2-6}
				&\multirow{2}*{160}&PDNR&44&1.43e+08&$\pm 10.3$\\
				&&PQNR&3&4.37e+07&$\pm 130.9$\\
				\midrule
				\multirow{4}*{\texttt{max\_outer\_iterations}}
				&\multirow{1}*{64}&PDNR&13&4.01e+07&$\pm 
				4.4$\\
				&\multirow{1}*{128}&PDNR&52&6.01e+07&$\pm 5.3$\\
				&\multirow{1}*{256}&PDNR&73&7.97e+07&$\pm 7.1$\\
				&\multirow{1}*{512}&PDNR&73&9.33e+07&$\pm 6.7$\\
				\midrule
				\multirow{12.5}*{\texttt{max\_backtrack\_steps}}
				&\multirow{2}*{1}&PDNR&1&4.82e+07&$\pm 0.0$\\
				&&PQNR&16&2.98e+08&$\pm 50.7$\\
				\cmidrule{2-6}
				&\multirow{2}*{2}&PDNR&17&6.68e+07&$\pm 16.7$\\
				&&PQNR&16&3.76e+08&$\pm 60.8$\\
				\cmidrule{2-6}
				&\multirow{2}*{4}&PDNR&17&8.01e+07&$\pm 17.6$\\
				&&PQNR&26&2.02e+08&$\pm 46.2$\\
				\cmidrule{2-6}
				&\multirow{2}*{8}&PDNR&20&8.93e+07&$\pm 16.0$\\
				&&PQNR&5&1.20e+08&$\pm 87.9$\\
				\cmidrule{2-6}
				&\multirow{2}*{10}&PDNR&74&8.94e+07&$\pm 6.8$\\
				&&PQNR&1&9.69e+07&$\pm 0.0$\\
				\cmidrule{2-6}
				&\multirow{1}*{16}&PDNR&15&2.20e+08&$\pm 18.9$\\
				\midrule
				\multirow{3}*{\texttt{min\_variable\_nonzero\_tolerance}}
				&\multirow{2}*{1e-07}&PDNR&65&9.08e+07&$\pm 7.1$\\
				&&PQNR&1&9.69e+07&$\pm 0.0$\\
				\cmidrule{2-6}
				&\multirow{1}*{0.001}
				&PDNR&21&9.02e+07&$\pm 12.0$\\
				\midrule
				\multirow{6}*{\texttt{size\_LBFGS}}
				&\multirow{1}*{3}&PQNR&1&9.69e+07&$\pm 0.0$\\
				&\multirow{1}*{4}&PQNR&1&9.12e+07&$\pm 0.0$\\
				&\multirow{1}*{5}&PQNR&1&1.73e+08&$\pm 0.0$\\
				&\multirow{1}*{10}&PQNR&1&4.68e+07&$\pm 0.0$\\
				&\multirow{1}*{15}&PQNR&1&1.41e+08&$\pm 0.0$\\
				&\multirow{1}*{20}&PQNR&2&6.06e+07&$\pm 27.1$\\
				\midrule
				\multirow{8}*{\texttt{step-reduction-factor}}
				&\multirow{1}*{0.1}&PDNR&27&1.53e+08&$\pm 17.4$\\
				\cmidrule{2-6}
				&\multirow{2}*{0.5}&PDNR&65&9.08e+07&$\pm 7.1$\\
				&&PQNR&1&9.69e+07&$\pm 0.0$\\
				\cmidrule{2-6}
				&\multirow{2}*{0.7}&PDNR&3&6.41e+07&$\pm 30.5$\\
				&&PQNR&13&2.21e+08&$\pm 45.3$\\
				\cmidrule{2-6}
				&\multirow{2}*{0.9}&PDNR&14&1.31e+08&$\pm 16.8$\\
				&&PQNR&17&3.77e+08&$\pm 46.9$\\
				\midrule
				\multirow{5}*{\texttt{suff\_decrease\_tolerance}}
				&\multirow{2}*{0.0001}&PDNR&65&9.08e+07&$\pm 7.1$\\
				&&PQNR&1&9.69e+07&$\pm 0.0$\\
				\cmidrule{2-6}
				&\multirow{2}*{0.01}&PDNR&20&9.35e+07&$\pm 11.7$\\
				&&PQNR&2&3.75e+07&$\pm 666.2$\\
				\bottomrule
	\end{tabular}}}%
\end{table}
\label{app:supplementary}
\clearpage

\end{document}